\documentclass[12pt]{article}
\usepackage{amsfonts}
\usepackage{amssymb}
\usepackage{amsmath,amsthm}
\usepackage{enumerate}

\usepackage{comment}

\usepackage{bbm}

\newtheorem{lemma}{Lemma}
\newtheorem{theorem}{Theorem}
\newtheorem{theorema}{Theorem} 

\newtheorem{lemmaa}{Lemma} 

\theoremstyle{definition}\newtheorem{Example}{Example}\theoremstyle{plain}

\theoremstyle{definition}
\newtheorem{assumption}{Assumption}
\newtheorem{assumptionA}{Assumption}
\theoremstyle{plain}

\renewcommand{\hat}{\widehat}
\renewcommand{\tilde}{\widetilde}

\newcommand{\ind}[1]{\mathbbm{1}\!\left\{#1\right\}}

\newcommand{\tttt}[1]{#1}

\input{tcilatex}

\begin{document}

\title{New goodness-of-fit diagnostics for conditional discrete response
models}
\author{\textsc{Igor Kheifets}\thanks{%
ITAM, Mexico. Email: igor.kheifets@itam.mx} \textsc{\ \ and } \textsc{Carlos
Velasco}\thanks{%
Department of Economics, Universidad Carlos III de Madrid. Email:
carlos.velasco@uc3m.es} }
\maketitle

\begin{abstract}
This paper proposes new specification tests for conditional models with
discrete responses, which are key to apply efficient maximum likelihood
methods, to obtain consistent estimates of partial effects and to get
appropriate predictions of the probability of future events. In particular,
we test the static and dynamic ordered choice model specifications and can
cover infinite support distributions for e.g. count data. The traditional
approach for specification testing of discrete response models is based on
probability integral transforms of a jittered discrete data which leads to
continuous uniform iid series under the true conditional distribution. Then,
standard specification testing techniques for continuous variables could be
applied to the transformed series, but the extra randomness from jitters
affects the power properties of these methods. We investigate in this paper
an alternative transformation based only on original discrete data that
avoids any randomization. We analyze the asymptotic properties of
goodness-of-fit tests based on this new transformation and explore the
properties in finite samples of a bootstrap algorithm to approximate the
critical values of test statistics which are model and parameter dependent.
We show analytically and in simulations that our approach dominates the
methods based on randomization in terms of power. We apply the new tests to
models of the monetary policy conducted by the Federal Reserve.\newline


\end{abstract}

\pagebreak

\renewcommand{\baselinestretch}{1.5}\normalsize

\section{INTRODUCTION}

Many statistical models specify the conditional distribution of a discrete
response variable given some explanatory variables, including the
description of binary, multinomial, ordered choice and count data. In this
paper we analyze goodness-of-fit tests for both static models with
covariates as well as dynamic ordered choice and count data models, where
the conditioning information set may also include past information on the
discrete variable and a set of (contemporaneous) explanatory variables which
frequently appear in the social sciences, see Kedem and Fokianos (2002) and
Greene and Hensher (2010). For example, dynamic models are popular in
macroeconomic applications, see for instance Hamilton and Jord\'{a} (2002),
Dolado and Maria-Dolores (2002) and Basu and de Jong (2007) for modeling
central banks decisions or Kauppi and Saikkonen (2008) and Startz (2008) for
predicting US recessions; in finance, see e.g. Rydberg and Shephard (2003)
for modeling the size of asset price movements and Fokianos et al.\ (2009) for
the number of transactions per minute of a particular stock.

Suppose we observe the random variables $\{Y_{t},X^{\prime }_{t}\}_{t=1}^T$
and consider the information sets $\Omega _{t}=\left \{ X_{t}, Y_{t-1},
X_{t-1},Y_{t-2},X_{t-2},\ldots \right \}$ for each period $t=1,2,\ldots,T$.
We are interested in testing the null hypothesis that the distribution of $%
Y_{t}$ conditional on $\Omega _{t}$ is in the parametric family $F_{t,\theta
}(\cdot \mid \Omega _{t})$, i.e.%
\begin{equation*}
H_{0}:Y_{t}\mid \Omega _{t}\sim F_{t,\theta _{0}}(\cdot \mid \Omega _{t})\ \ 
\text{for\ some }\theta _{0}\in \Theta ,\ t=1,2,\ldots,T,
\end{equation*}%
where $\Theta \subset \mathop{\mathbb R}\nolimits^{m}$ is the parameter
space, while the alternative hypothesis ($H_1$) for the omnibus test would
be the negation of $H_{0}$.

We consider a class $\mathop{\mathcal M}\nolimits$ of discrete conditional
distributions defined on $\mathop{\mathcal K}\nolimits=\{1,2,\ldots ,K\}$,
for integer $K>1$ or on $\mathop{\mathcal K}\nolimits=\{1,2,\ldots ,\infty
\} $ such that for all $F\in \mathop{\mathcal M}\nolimits$ it holds that $%
F\left( 0\right) =0$, $f\left( k\right) :=F\left( k\right) -F\left(
k-1\right) >0$ for all $k=1,2,\ldots $ and $\sum_{k\in \mathop{\mathcal K}%
\nolimits}f(k)=1$. 
This setup includes numerous models that have been used extensively in
applied work both for dynamic and for iid data, here we describe briefly two
of them.

\bigskip

\begin{Example}[Dynamic multinomial ordered choice model]
The discrete responses $Y_{t}$ are assumed to be generated by the rule%
\begin{equation*}
Y_{t}=\left \{ 
\begin{array}{ccc}
1 & \ \text{if} & V_{t}^{\ast }\leq \tau _{1} \\ 
2 & \ \text{if} & \tau _{1}<V_{t}^{\ast }\leq \tau _{2} \\ 
& \vdots &  \\ 
K & \ \text{if} & V_{t}^{\ast }>\tau _{K-1},%
\end{array}%
\right.
\end{equation*}%
where $V_{t}^{\ast }$ is a continuous latent variable and $\tau _{1},\ldots
,\tau _{K-1}$ are threshold parameters that define $K$ intervals in $\mathbb{%
R}$. In a simple model, e.g. Basu and de Jong (2007), the latent variable is
determined through the linear equation%
\begin{equation*}
V_{t}^{\ast }=X_{t}^{\prime }\beta +\rho Y_{t-1}+\varepsilon _{t},
\end{equation*}%
where $X_{t}$ is a vector of stationary exogenous regressors, $\beta $ a
vector of regression parameters, $\varepsilon _{t}$ is the shock in each
period, and $Y_{t-1}$ could be replaced by any function of the past $\left
\{ Y_{t-1},\ldots ,Y_{t-n}\right \} $ for some finite $n.$ The cdf of $%
\varepsilon _{t},$ $F_{\varepsilon },$ is going to determine the class of
multinomial model, i.e. ordered multinomial probit (if $\varepsilon _{t}$ is
standard normal) or logit (if $\varepsilon _{t}$ is logistic), since $%
F_{t,\theta _{0}}$ is defined at once from%
\begin{equation*}
\begin{split}
\Pr \left( Y_{t}=k\mid \Omega _{t}\right) & =\Pr \left( \tau
_{k-1}<V_{t}^{\ast }\leq \tau _{k}\mid \Omega _{t}\right) \\
& =F_{\varepsilon }\left( \tau _{k}-X_{t}^{\prime }\beta -\rho
Y_{t-1}\right) -F_{\varepsilon }\left( \tau _{k-1}-X_{t}^{\prime }\beta
-\rho Y_{t-1}\right) ,
\end{split}%
\end{equation*}%
with $\tau _{0}=-\infty $ and $\tau _{K}=\infty \ $and $\theta _{0}=\left(
\beta ^{\prime },\rho ,\tau _{1},\ldots ,\tau _{K-1}\right) ^{\prime }.$
\end{Example}

\bigskip

\begin{Example}[Poisson Model]
The variate $Y_{t}=Y_{t}^{\ast }+1$ is defined on the counts $Y_{t}^{\ast
}=0,1,2,\ldots $ which are assumed to follow a conditional Poisson
distribution%
\begin{equation*}
Y_{t}^{\ast }\mid \Omega _{t}\sim \text{Poisson}(\lambda _{t}),
\end{equation*}%
where the conditional mean can depend on covariates through an exponential link as $\lambda _{t}=\exp(X_t'\beta)$ or on previous observations
through an identity link as $\lambda _{t}=\alpha _{0}+\alpha _{1}\lambda
_{t-1}+\rho Y_{t-1}^{\ast },$ e.g. Fokianos et al.\ (2009), or through the
logarithmic canonical link as log$\left( \lambda _{t}\right) =X_{t}^{\prime
}\beta +\rho e_{t-1},$ where $e_{t}=\left( Y_{t}^{\ast }-\lambda _{t}\right)
/\lambda _{t}$ are scaled and centered errors, e.g. Davis et al.\ (2003).
\end{Example}

\bigskip

Despite that a correct specification is key to apply efficient maximum
likelihood methods, to obtain consistent estimates of partial effects and to
get appropriate predictions of the probability of future events, empirical
researchers typically do not perform goodness of fit testing of such models
as they would do in a continuous case. In general, there are only a few
specification tests available for discrete data, see Mora and Moro-Egido
(2007). Two of them, the test of the Generalized Linear Model (GLM) of Stute
and Zhu (2002) and the conditional Kolmogorov test of Andrews (1997), based
on the specification of the conditional mean for binary data, can be adapted
for this purpose and we discuss this possibility and compare it to our
approach in Section~6. A related test to Andrews derived for time
series by Corradi and Swanson (2006) could be adapted also for discrete
data, but this is testing a different null hypothesis concerning a
distribution given a finite conditioning set not characterizing the complete
dynamics of the process. There are also tests designed specifically for Poisson models (see e.g.\ Neumann 2011; Fokianos and Neumann, 2013).

In what follows we propose conditional, dynamic discrete analogs of
the Kolmogorov-Smirnov goodness of fit measure that can exploit different
restrictions derived from the martingale difference property of a particular
transformation of the data under the null hypothesis. This property is
derived from the specification of a complete dynamic model given the
information set generated by all the past observations of the discrete
response and other explanatory variables and is used to build the asymptotic
theory for our tests. Under i.i.d.\ assumptions this martingale
difference property leads to an exact independence of the transformation
sequence under the null and a much simpler parallel asymptotic theory.

When the fitted distribution is continuous, the \emph{relative distribution}
of $Y_{t}$ compared to $F_{t,\theta _{0}}$ defined as the cdf of the
Rosenblatt's (1952) transforms, also called conditional Probability Integral
Transforms (PIT),%
\begin{equation*}
U_{t}\left( \theta _{0}\right) :=F_{t,\theta _{0}}\left( Y_{t}\mid \Omega
_{t}\right) ,\qquad t=1,2,\ldots ,T
\end{equation*}%
is standard uniform and $U_{t}\left( \theta _{0}\right)$ are distributed as independent $[0,1]$ uniform random variables under $H_{0}$%
. This serves as a basis for several specification tests of $H_{0}$, see
e.g. Bai (2003) and Kheifets (2015) for dynamic models and Delgado and Stute
(2008) for independent and identical distributed (iid) data. However
Rosenblatt transformation is not appropriate for discrete support random
variables, producing non-iid pseudo residuals even under the null of correct
specification. To solve the limitations of PIT-based testing techniques for
discrete data, several alternative transforms have been proposed, see Jung,
Kukuk and Liesenfeld (2006), Czado, Gneiting and Held (2009) and references
therein. An easy and popular way is to randomize, i.e. to interpolate the
discrete values of $Y_{t}$ with independent noise in $[0,1]$, recent
references include Kheifets and Velasco (2013) and Lee (2014). Unfortunately
the additional simulated noise affects the power of the tests and may lead
to different conclusions depending on the simulation outcome.

In this paper instead, we consider a \emph{nonrandomized transform} $%
Y_{t}\mapsto I_{t,\theta _{0}}\left( u\right) $ for $u\in \lbrack 0,1]$, 
\begin{equation}
I_{t,\theta _{0}}\left( u\right) :=\left\{ 
\begin{array}{rrr}
0, &  & u\leq U_{t}^{-}\left( \theta _{0}\right) ; \\ 
\displaystyle{\ \frac{u-U_{t}^{-}\left( \theta _{0}\right) }{U_{t}\left(
\theta _{0}\right) -U_{t}^{-}\left( \theta _{0}\right) },} &  & 
U_{t}^{-}\left( \theta _{0}\right) \leq u\leq U_{t}\left( \theta _{0}\right)
; \\ 
1, &  & U_{t}\left( \theta _{0}\right) \leq u,%
\end{array}%
\right.   \label{eq:nr}
\end{equation}%
where $U_{t}^{-}\left( \theta _{0}\right) :=F_{t,\theta _{0}}\left(
Y_{t}-1\mid \Omega _{t}\right) $. This transform, conditional on data, is
nonrandomized in the sense that it does not depend on extra sources of
randomness, as opposed to interpolation transforms discussed in the next
section. The unconditional version of this transform appears in Handcock and
Morris (1999) and more recently in Czado, Gneiting and Held (2009)
where it is used for calibration, but no formal tests are proposed
there. This transformation can also be seem as a particular case of the multilinear extension as defined in Genest, Ne\v slehov\' a and R\' emillard (2014).
 As we show below, for every $u\in \lbrack 0,1]$, $I_{t,\theta _{0}}\left( u\right) -u$
constitute a martingale difference sequence (MDS) with respect to $\Omega
_{t}$ under $H_{0}$ and can be used for testing $H_{0}$ as $I_{t,\theta
_{0}}\left( u\right) $ loses this property when the model is misspecified.
For instance, we can compute the pseudo empirical relative distribution of $%
Y_{t}$ compared to $F_{t,\theta _{0}}$
\begin{equation*}
\tilde{F}_{\theta _{0}}\left( u\right) :=\frac{1}{T}\sum_{t=1}^{T}I_{t,%
\theta _{0}}\left( u\right) ,\ \ \ u\in \left[ 0,1\right] ,
\end{equation*}%
which can be contrasted with the uniform cdf using the following empirical
process 
\begin{equation*}
S_{1T}\left( u\right) :=\frac{1}{T^{1/2}}\sum_{t=1}^{T}\left\{ I_{t,\theta
_{0}}\left( u\right) -u\right\} ={T^{1/2}}\left( \tilde{F}_{\theta
_{0}}\left( u\right) -u\right) ,
\end{equation*}%
which converges weakly to a Gaussian process. In addition, in order to
control dynamics in $I_{t,\theta _{0}}\left( u\right) $, we can compare the
joint pseudo empirical cdf with the uniform on a square using the
biparameter process 
\begin{equation}
S_{2T}\left( u\right) :=\frac{1}{(T-1)^{1/2}}\sum_{t=2}^{T}\left\{
I_{t,\theta _{0}}\left( u_{1}\right) I_{t-1,\theta _{0}}\left( u_{2}\right)
-u_{1}u_{2}\right\} ,  \label{eq:S2}
\end{equation}%
where $u=\left( u_{1},u_{2}\right) $. To obtain feasible tests we need to
consider norms of $S_{jT}$ for $j=1,2$. We use the Cramer-von Mises $\int
S_{jT}\left( u\right) ^{2}d\varphi \left( u\right) $ for some absolute
continuous measure~$\varphi $ in $\left[ 0,1\right] ^{j}$, or
Kolmogorov-Smirnov $\sup_{u\in \lbrack 0,1]^{j}}\left\vert S_{jT}\left(
u\right) \right\vert $ norms.

When the parameter $\theta _{0}$ is unknown under the null, we use an
estimate $\hat{\theta}_{T}$ and account for the parameter estimation effect
in the $p$-value computation with a parametric bootstrap method. It might be
possible also to derive, e.g. martingale, distribution-free transforms, but
since they typically need to be programmed on a case by case basis for each
model, so can be impractical, and are beyond the scope of this paper. As
far as we know, our proposal is the first formal specification test of
ordered discrete choice models which accounts properly for parameter
uncertainty and is based on a nonrandomized transform, which makes it
attractive in terms of power against a wide set of alternative hypotheses.

The rest of the paper is organized as follows. In the next section, we
describe different alternatives to the PIT. In Sections 3 and 4, we provide
the main asymptotic properties of the nonrandomized transforms and of the
resulting univariate and bivariate empirical processes using martingale
theory. In particular, we establish weak limits under fixed and local
alternatives accounting for parameter estimation effect. Section 5 discusses
the implementation of new tests with a simple bootstrap algorithm. Section 6
provides a small simulation exercise and an application exploring the
properties of specification tests based on both randomized and non
randomized transformations. Then we conclude. All proofs are contained in
the Appendix.

\section{ALTERNATIVES TO PIT FOR DISCRETE DATA}

In order to further motivate the nonrandomized transform $I_{t,\theta _{0}}$
defined in~(\ref{eq:nr}), we introduce the \textit{randomized} PIT, 
\begin{equation}
U_{t}^{r}\left( \theta _{0}\right) :=U_{t}^{-}\left( \theta _{0}\right)
+Z_{t}^{U}\left( U_{t}\left( \theta _{0}\right) -U_{t}^{-}\left( \theta
_{0}\right) \right) ,  \label{eq:Ur}
\end{equation}%
where $\{Z_{t}^{U}\}_{t=1}^{T}$ are independent standard uniform random
variables, and independent of $Y_{t}$. Alternatively, $U_{t}^{r}$ can be
obtained by applying the standard continuous PIT to the \textit{continuous}
random variable $Y_{t}^{\dag }:=Y_{t}-1+Z_{t}$, where $\{Z_{t}\}_{t=1}^{T}$
are iid with any continuous cdf $F_{Z}$ on $[0,1]$. Indeed, we can construct
the cdf of $Y_{t}^{\dag }$, 
\begin{equation*}
F_{t,\theta _{0}}^{\dag }\left( y\mid \Omega _{t}\right) =F_{t,\theta
_{0}}(\left\lfloor y\right\rfloor \mid \Omega _{t})+F_{Z}\left(
y-\left\lfloor y\right\rfloor \right) \left( F_{t,\theta _{0}}\left(
\left\lfloor y+1\right\rfloor \mid \Omega _{t}\right) -F_{t,\theta
_{0}}\left( \left\lfloor y\right\rfloor \mid \Omega _{t}\right) \right) ,
\end{equation*}%
where $\left\lfloor y\right\rfloor $ is the floor function, i.e. the maximum
integer not exceeding $y$, and find that 
\begin{equation*}
U_{t}^{r}\left( \theta _{0}\right) =F_{t,\theta _{0}}^{\dag }\left(
Y_{t}^{\dag }\mid \Omega _{t}\right) ,
\end{equation*}%
for $Z_{t}^{U}=F_{Z}\left( Z_{t}\right) $ and any choice of $F_{Z}$, see
Kheifets and Velasco (2013). Note that the cdf of $Y_{t}^{\dag }$
conditional on $\Omega _{t}$ and $\left\{ \Omega _{t},Z_{t-1},Z_{t-2},\ldots
,Z_{1}\right\} $ coincide. Under $H_{0}$, $U_{t}^{r}\left( \theta
_{0}\right) $ are iid $U\left[ 0,1\right] $ variables as under any
continuous distribution specification, while $U_{t}\left( \theta _{0}\right) 
$ and $U_{t}^{-}\left( \theta _{0}\right) $ are not independent nor $U\left[
0,1\right] $. Using the typical discrepancy measures, the empirical cdf of 
$U_{t}^{r}\left( {\theta }_{0}\right) $, estimated using the \emph{%
randomized transform} $Y_{t}\mapsto \mathbbm{1}\!\left\{ U_{t}^{r}\left(
\theta _{0}\right) \leq u\right\} $, 
\begin{equation*}
\hat{F}_{\theta _{0}}^{r}\left( u\right) :=\frac{1}{T}\sum_{t=1}^{T}%
\mathbbm{1}\!\left\{ U_{t}^{r}\left( \theta _{0}\right) \leq u\right\} ,\ \
\ u\in \left[ 0,1\right] ,
\end{equation*}%
can be compared to the uniform cdf. Kheifets and Velasco (2013) then test $%
H_{0}$ using empirical process based on the randomized transform 
\begin{equation*}
{R}_{1T}\left( u\right) :=T^{1/2}\left\{ \hat{F}_{{\theta }_{0}}^{r}\left(
u\right) -u\right\} =\frac{1}{T^{1/2}}\sum_{t=1}^{T}\left[ \mathbbm{1}%
\!\left\{ U_{t}^{r}\left( {\theta }_{0}\right) \leq u\right\} -u\right] ,\ \
\ u\in \left[ 0,1\right] .
\end{equation*}

We can also consider reducing the dependence on a particular outcome of the
noise $Z_{t}^{U}$ in (\ref{eq:Ur}) and in the randomized transform by taking
averages over $M$ replications of $\{Z_{t}^{U}\}_{t=1}^{T}$, conditional on
the original data, similar to \textquotedblleft average-jittering" of
Machado and Santos Silva (2005). Suppose that for each $t$ we have $M$
independent sequences of uniform $U[0,1]$ noises $Z_{t,m}^{U}$, $%
m=1,2,\ldots ,M$, which generate $U_{t,m}^{r}\left( \theta _{0}\right) $
according to (\ref{eq:Ur}). Define the \textit{M-random} transform $%
Y_{t}\mapsto I_{t,\theta _{0},M}\left( Y_{t},u\right) $, 
\begin{equation*}
I_{t,\theta _{0},M}\left( Y_{t},u\right) :=\frac{1}{M}\sum_{m=1}^{M}%
\mathbbm{1}\!\left\{ U_{t,m}^{r}\left( {\theta }_{0}\right) \leq u\right\} ,
\end{equation*}%
which takes values on the set $\{0,1/M,2/M,\ldots ,1\}$ and has mean $u$
under $H_{0}$. Then the cdf of $U_{t}^{r}\left( {\theta }_{0}\right) $ is
estimated by 
\begin{equation*}
\hat{F}_{\theta _{0},M}^{r}\left( u\right) :=\frac{1}{T}\sum_{t=1}^{T}I_{t,%
\theta _{0},M}\left( Y_{t},u\right) ,\ \ \ u\in \left[ 0,1\right] .
\end{equation*}%
Note that with $M=1$ we are back to $\hat{F}_{\theta _{0}}^{r}\left(
u\right) $, and therefore, we can generalize ${R}_{1T}$ to 
\begin{equation*}
{R}_{1T,M}\left( u\right) :=T^{1/2}\left\{ \hat{F}_{{\theta }%
_{0},M}^{r}\left( u\right) -u\right\} ,\ \ \ u\in \left[ 0,1\right] .
\end{equation*}

In order to propose specification tests, following Handcock and Morris
(1999), we define the \emph{discrete relative distribution} of $Y_{t}$
compared to $F_{t,\theta _{0}}$ as the cdf of $U_{t}^{r}\left( {\theta }%
_{0}\right) $. Under $H_{0}$, the discrete relative distribution is the
uniform $U\left[ 0,1\right] $. As we show in the next section, three
consistent estimators of the discrete relative distribution of $Y_{t}$
compared to $F_{t,\theta _{0}}$ can be ordered in terms of efficiency in the
following way: $\tilde{F}_{\theta _{0}}\left( u\right) $ (the most
efficient), $\hat{F}_{\theta _{0},M}^{r}\left( u\right) $ and $\hat{F}%
_{\theta _{0}}^{r}\left( u\right) $. This ordering is determined by the
amount of noise introduced in the definitions of the transforms: i.e. in
nonrandomized, $M$-randomized and ($1$-)randomized transforms. The
nonrandomized transform can be equivalently obtained by integrating out the
extra noise in the randomized transform $I_{t,\theta _{0}}\left(
Y_{t},u\right) =\int \mathbbm{1}\!\left\{ U_{t}^{r}\left( \theta _{0}\right)
\leq u\right\} dF_{Z}$ or by taking the number of replications $M$ to
infinity, thus completely removing the noise from the estimate of the
discrete relative distribution and other functionals of the transforms. The
efficiency of the nonrandomized transform translates into the increased
power of the specification tests based on this transform, whose properties
we study next.

\section{PROPERTIES OF EMPIRICAL PROCESSES BASED ON THE NONRANDOMIZED
TRANSFORM}

As shown in the next lemma, the building blocks of $\tilde{F}_{\theta
_{0}}\left( u\right) ,$ $I_{t,\theta _{0}}\left( u\right) -u$, constitute a
martingale difference sequence (MDS) with respect to $\Omega _{t}$, and
therefore $\tilde{F}_{\theta _{0}}\left( u\right) $ is an unbiased and
consistent estimate of the uniform cdf under the null, a reasonable basis
for developing tests of $H_{0}$. Moreover, the MDS property will allow us to
establish the asymptotic properties of our test without imposing any
additional restrictions. Let for $u,v\in \lbrack 0,1]$ 
\begin{equation*}
\gamma _{t,\theta _{0}}\left( u,v\right) :=\frac{\left( F_{k}-u\vee v\right)
\left( u\wedge v-F_{k-1}\right) }{F_{k}-F_{k-1}}\mathbbm{1}\!\left\{
F_{t,\theta _{0}}^{-1}\left( u\mid \Omega _{t}\right) =F_{t,\theta
_{0}}^{-1}\left( v\mid \Omega _{t}\right) \right\} ,
\end{equation*}%
where $k=k\left( u\right) =F_{t,\theta _{0}}^{-1}\left( u\mid \Omega
_{t}\right) $, with $F_{t,\theta _{0}}^{-1}\left( u\mid \Omega _{t}\right)
:=\min \{y:F_{t,\theta _{0}}\left( y\mid \Omega _{t}\right) \geq u\}$ being
the conditional quantile function and $F_{k}:=F_{t,\theta _{0}}\left( k\mid
\Omega _{t}\right) $.

\begin{lemma}
\label{lem:mds} Under $H_{0}$, $I_{t,\theta _{0}}\left( u\right) -u$ is a
martingale difference sequence with respect to $\Omega _{t}$, i.e.%
\begin{equation*}
\mathop{\rm E}\nolimits\left[ I_{t,\theta _{0}}\left( u\right) \mid \Omega
_{t}\right] =u,\ \ \ a.s.,
\end{equation*}%
%
%
%
%
%
%
%
%
%
%
%
with conditional covariance 
\begin{equation*}
\mathop{\rm E}\nolimits\left[ I_{t,\theta _{0}}\left( u\right) I_{t,\theta
_{0}}\left( v\right) \mid \Omega _{t}\right] =u\wedge v-uv-\gamma _{t,\theta
_{0}}\left( u,v\right) ,\ \ \ a.s.
\end{equation*}
\end{lemma}

Note that $I_{t,\theta _{0}}\left( u\right) $ are not necessarily
independent across $t$ despite the fact that by the martingale difference
property, $I_{t,\theta _{0}}\left( u\right) $ and $I_{t-j,\theta _{0}}\left(
v\right) $ are serially uncorrelated for all $j\neq 0$ and all $u,v\in \left[
0,1\right] ,$ see the Appendix. On the other hand, the $I_{t,\theta
_{0}}\left( u\right) $ are (conditionally) heteroskedastic, therefore the
variance of $S_{1T}$ is model and parameter dependent, but its distribution
can be simulated conditional on exogenous information in $\Omega _{t}.$ 

Let $V_{1T}\left( u,v\right) :=\mathop{\rm Cov}\nolimits\left[ S_{1T}\left(
u\right) ,S_{1T}\left( v\right) \right] $, then since $0\leq \gamma
_{t,\theta _{0}}\left( u,v\right) <1\ a.s.$, 
\begin{equation*}
V_{1T}\left( u,v\right) =u\wedge v-uv-\frac{1}{T}\sum_{t=1}^{T}\mathop{\rm E}%
\nolimits\left[ \gamma _{t,\theta _{0}}\left( u,v\right) \right] \leq
u\wedge v-uv,
\end{equation*}%
i.e. the covariance and variance of ${S}_{1T}$ are not larger than those of
the randomized transformation-based process ${R}_{1T}$ or its weak limit,
the Brownian sheet, see Corollary 4 in Kheifets and Velasco (2013). 

Due to Lemma \ref{lem:mds}, $\mathop{\rm E}\nolimits \left[ \tilde{F}%
_{\theta _{0}}\left( u\right) \right] =u$ under $H_{0}$ and the natural
empirical process for performing tests on $H_{0}$ is then $S_{1T}$. This
process, being based on a nonrandomized transform, does not involve the
extra noise that appears in the randomized transform based empirical process 
${R}_{1T}$ for testing $U_{t}^{r}\sim U[0,1]$, proposed by Kheifets and
Velasco (2013), or in its modification ${R}_{1T,M}$ based on the $M$%
-randomized transform. The next lemma is the key to understand the
improvement of the $M$-randomized over the randomized and of the
nonrandomized, advocated in this paper, over the $M$-randomized transform
approaches.

\begin{lemma}
\label{lem:comp} Suppose that the uniform law of large numbers holds for $%
\hat{F}_{\theta _{0},M}^r\left( u\right)$ and $\tilde{F}_{\theta _{0}}\left(
u\right)$. Independently of whether $H_0$ holds or not, $\hat{F}_{\theta
_{0},M}^r\left( u\right)$ and $\tilde{F}_{\theta _{0}}\left( u\right)$
consistently and uniformly in $u$ estimate the relative distribution, i.e.
the cdf of $U^r_{t}\left( \theta _{0}\right)$. $\tilde{F}_{\theta
_{0}}\left( u\right)$ is more efficient, but the difference in efficiency
goes to $0$ as $M\to \infty$. In particular, under $H_0$, 
\begin{equation*}
\mathop{\rm E}\nolimits \left[R_{1T,M}\left(u\right)R_{1T,M}\left(v\right)%
\right]=\frac{1}{M} \mathop{\rm E}\nolimits \left[ R_{1T}\left(u%
\right)R_{1T}\left(v\right)\right]+\left(1-\frac{1}{M}\right) \mathop{\rm E}%
\nolimits\left[S_{1T}\left(u\right)S_{1T}\left(v\right)\right].
\end{equation*}
\end{lemma}

From Lemma~\ref{lem:comp}, it follows that $S_{1T}$ has the smallest
variance, the variance of $R_{1T,M}$ is a weighted sum of those of $S_{1T}$
and $R_{1T}$, see also Equation (5) in Machado and Santos Silva (2005).
Other advantages of $S_{1T}$ over $R_{1T,M}$, are 1) computational, as there
is no need to simulate $M$ paths of transformations and 2) theoretical,
since the weak convergence is easier to prove for processes which are
piece-wise linear in parameters. Therefore we concentrate on studying the
properties of tests based on the nonrandomized transform, for which we
introduce the following assumption.

\begin{assumption}
${F_{t,\theta _{0}}\left( \cdot \mid \Omega _{t}\right) }%
\left( k\right) \in \mathop{\mathcal M}\nolimits$ $a.s.$ for all $t$.
Moreover, there exists a finite function $\gamma _{\infty }\left( u,v\right) 
$, such that uniformly in $\left( u,v\right) \in \left[ 0,1\right] ^{2},$ $%
T^{-1}\sum_{t=1}^{T}\gamma _{t,\theta _{0}}\left( u,v\right) \rightarrow
_{p}\gamma _{\infty }\left( u,v\right) $.
\end{assumption}

This assumption implicitly restrict dynamics such that a uniform law
of large numbers (LLN) holds for the averaged conditional covariance
function. In the case of stationary and ergodic data, $\gamma _{\infty
}\left( u,v\right) =\mathop{\rm E}\nolimits \left[ \gamma _{1,\theta
_{0}}\left( u,v\right) \right]$. Sufficient conditions for the
stationarity and ergodicity of dynamic multinomial ordered choice models are
given in Basu and de Jong (2007) and for autoregressive Poisson are given in
Davis et al.\ (2003), Fokianos et al.\ (2009) and Doukhan et al.\ (2012). Then
it is possible to show the uniformity of the convergence from a point-wise
result, since the summands are continuous, piece-wise polynomials in
$u$ and $v$. As an illustration, in Section 8.5 in  Appendix we discuss the assumptions for the Poisson model.

The next result describes the asymptotic distribution of $S_{1T}$ under the
null hypothesis. Let $\Rightarrow $ denote weak convergence in $\ell
^{\infty }\left[ 0,1\right] $, see e.g. van der Vaart and Wellner (1996). In
fact, our empirical processes are continuous, which simplifies tightness
verification. Let $V_{1\infty }\left( u,v\right) :=u\wedge v-uv-\gamma
_{\infty }\left( u,v\right) $.

\begin{lemma}
\label{lem:limitS} Suppose Assumption 1 holds. Under $H_{0}$,%
\begin{equation*}
S_{1T}\Rightarrow S_{1\infty},
\end{equation*}%
where $S_{1\infty}$ is a Gaussian process in $\left[ 0,1\right] $ with zero
mean and covariance function $V_{1\infty}$.
\end{lemma}

The asymptotic distribution of $S_{1 T}$ is model and parameter dependent,
and the practical implementation of tests when $\theta _{0}$ is unknown is
discussed in Section 3.2 after presenting a general class of local
alternatives to the null of correct specification of the conditional
distribution.

\subsection{Local Alternatives}

We next discuss the asymptotic properties of the empirical process $S_{1T}$
under a class of alternative hypothesis, that will lead to consistency of
the specification tests based on $S_{1T}$ for a wide class of alternatives.
We consider the following class of local alternatives to $H_{0},$%
\begin{equation*}
H_{1T}:\ \ Y_{t}\mid \Omega _{t}\sim G_{T,t,\theta _{0}}(\cdot \mid \Omega
_{t})\ \ \text{for some }\theta _{0}\in \Theta ,
\end{equation*}%
where 
\begin{equation*}
G_{T,t,\theta _{0}}(y\mid \Omega _{t})=\left( 1-\frac{\delta }{T^{1/2}}%
\right) F_{t,\theta _{0}}(y\mid \Omega _{t})+\frac{\delta }{T^{1/2}}%
H_{t}(y\mid \Omega _{t}),
\end{equation*}%
for some $0<\delta <T^{1/2}$ and for all $t$, $H_{t}\left( \cdot \mid \Omega
_{t}\right) \in \mathop{\mathcal M}\nolimits$. When $\delta =0$ then $H_{1T}$
nests $H_{0}.$

Following Kheifets and Velasco (2013), for any discrete distributions $G$
and $F$ in $\mathop{\mathcal M}\nolimits$, with probability functions $g$
and $f$, define 
\begin{equation*}
\begin{split}
d\left( G,F,u\right) ={}& G\left( F^{-1}\left( u\right) \right) -F\left(
F^{-1}\left( u\right) \right) \\
& -\frac{F\left( F^{-1}(u)\right) -u}{f\left( F^{-1}(u)\right) }\left[
g\left( F^{-1}(u)\right) -f\left( F^{-1}(u)\right) \right] .
\end{split}%
\end{equation*}%
Note, that $d\left( G,F,u\right) =E_G[I_F(Y,u)]-E_F[I_F(Y,u)]=E_G[I_F(Y,u)]-u$ and $d\left( G,F,u\right) \equiv 0$ if and only if $G\equiv F$. Under
any $G_{t}\left( \cdot \mid \Omega _{t}\right) \in \mathop{\mathcal M}%
\nolimits$, 
\begin{equation*}
\frac{1}{T^{1/2}}\mathop{\rm E}\nolimits \left[ S_{1T}\left( u\right) \right]
=\frac{1}{T}\sum_{t=1}^{T}\mathop{\rm E}\nolimits \left[ d\left( G_{t}\left(
\cdot \mid \Omega _{t}\right) ,F_{t,\theta _{0}}\left( \cdot \mid \Omega
_{t}\right) ,u\right) \right] .
\end{equation*}%
The next assumption guarantees that a LLN can be applied to the
empirical discrepancy between $H_{t}$ and $F_{t,\theta _{0}}.$

\begin{assumption}
Under $H_{1T}$, there exists a finite function $D_{1}\left( u\right) $, such
that uniformly in $u\in \left[ 0,1\right] ,\ \frac{1}{T}\sum_{t=1}^{T}d%
\left( H_{t}\left( \cdot \mid \Omega _{t}\right) ,F_{t,\theta _{0}}\left(
\cdot \mid \Omega _{t}\right) ,u\right) \rightarrow _{p}D_{1}\left( u\right) 
$.
\end{assumption}

Then the following lemma shows that the departure of $H_0$ in the direction
of $H_{1T}$ introduces a drift in the asymptotic distribution of $S_{1T}$
that will render consistency of hypothesis tests based on functionals of $%
H_{1T}$.

\begin{lemma}
\label{lem:limitSLA} Suppose Assumptions 1-2 hold. Under $H_{1T}$,%
\begin{equation*}
S_{1T}\Rightarrow S_{1\infty} +\delta D_1,
\end{equation*}%
where $S_{1\infty}$ is as in Lemma \ref{lem:limitS}.
\end{lemma}

\subsection{Parameter Estimation Effect}

In practice, tests based on $S_{1T}$ are unfeasible since $\theta _{0}$ is
unknown, and has to be estimated by $\hat{\theta}_{T},$ say. We assume that
we have available an estimate $\hat{\theta}_{T}$ so that under $H_{1T}$%
\begin{equation*}
T^{1/2}\left( \hat{\theta}_{T}-\theta _{0}\right) =O_{p}\left( 1\right) ,
\end{equation*}%
and define the process with estimated parameters 
\begin{equation*}
\hat{S}_{1T}\left( u\right) :=\frac{1}{T^{1/2}}\sum_{t=1}^{T}\left\{ I_{t,%
\hat{\theta}_{T}}\left( u\right) -u\right\} .
\end{equation*}%
We next analyze the consequences of replacing $\theta _{0}\ $by $\hat{\theta}%
_{T}$ in $\hat{S}_{1T}$.

Let $\Vert \cdot \Vert $ be Euclidean norm, i.e. for matrix $A$, $\Vert
A\Vert =\sqrt{\mathop{\rm tr}\nolimits \left( AA^{\prime }\right) }$, where $%
A^{\prime }$ is a transpose of $A$. For $\varepsilon >0,$ $B(a,\varepsilon )$%
\ is an open ball in $\mathop{\mathbb R}\nolimits^{m}$ with the center at
point $a$ and radius $\varepsilon $. For a cdf $F_{\theta }$ in $%
\mathop{\mathcal M}\nolimits$ define 
\begin{equation*}
\nabla \left( F_{\theta },u\right) :=\dot{F}_{{\theta }}\left( F_{\theta
}^{-1}\left( u\right) \right) -\frac{F_{{\theta }}\left( F_{\theta
}^{-1}\left( u\right) \right) -u}{{f_{\theta }}\left( F_{\theta }^{-1}\left(
u\right) \right) }\dot{f}_{\theta }\left( F_{\theta }^{-1}\left( u\right)
\right) ,
\end{equation*}%
where $\dot{F}_{\theta }:=\left( \partial /\partial \theta \right) F_{\theta
}$ and $\dot{f}_{\theta }:=\left( \partial /\partial \theta \right)
f_{\theta }$. We need the following assumptions to analyze the asymptotic
properties of $\hat{S}_{1T}$.

\begin{assumption}[Parametric family]
\begin{enumerate}
\item[(A)] The parameter space $\Theta $ is a compact set in a
finite-dimensional Euclidean space, $\theta \in \Theta \subset %
\mathop{\mathbb R}\nolimits^{m}$.

\item[(B)] There exists $\delta>0$, such that ${F_{t,\theta}}\left( \cdot
\mid \Omega _{t}\right)\in \mathop{\mathcal M}\nolimits$, for all $t$, $%
\Omega _{t}$, $T$ and $\theta \in B(\theta_0,\delta )$.

\item[(C)] $F_{t,\theta }\left( k\mid \Omega _{t}\right) $ is differentiable
with respect to $\theta \in B(\theta_0,\delta )$ and under $H_{1T}$ \newline
$\max_{t}\mathop{\rm E}\nolimits \left[ \max_{k}\sup_{\theta \in B(\theta_0,\delta )}\left \Vert 
\dot{F}_{t,\theta }\left( k\mid \Omega _{t}\right) \right \Vert \right] \leq
M_{F}<\infty .$

\item[(D)] Under $H_{1T}$, there exists a finite $L_{1}\left( u\right) :=%
\mathop{\rm plim}_{T\rightarrow \infty }T^{-1}\sum_{t=1}^{T}\nabla \left(
F_{t,\theta _{0}}\left( \cdot \mid \Omega _{t}\right) ,u\right) $.
\end{enumerate}
\end{assumption}

Conditions (A)-(C) about the parametric family of distribution are
standard, see e.g. Bai (2003, Assumptions A1-A2). For dynamic ordered choice
and Poisson models the differentiability of the conditional distribution
with respect to the parameter is equivalent to the differentiability of the
link function. Part (D) guarantees a nice limit behaviour of the average
generalized derivative of $I_{t,\theta }$. Conditions for no
effect of information truncation can be provided in a similar way to Bai
(2003, Assumption A4).

The following lemma provides an expansion of the empirical process with
estimated parameters as the sum of the process with known parameters and a
random drift describing parameter estimation.

\begin{lemma}
\label{lem:puexpansion} Suppose Assumptions 1-3 hold and $T^{1/2}\left(\hat{%
\theta}_{T}-\theta _{0}\right)=O_p(1)$. Under $H_{1T}$,%
\begin{equation}  \label{eq:puexpansion}
\hat{S}_{1T}\left( u\right) =S_{1T}\left( u\right) +T^{1/2}\left( \hat{\theta%
}_{T}-\theta _{0}\right) ^{\prime }\frac{1}{T}\sum_{t=1}^{T} \nabla
\left(F_{t,\theta_0}\left(\cdot \mid \Omega_t\right),u\right) +o_{p}\left(
1\right) ,
\end{equation}%
uniformly in $u$.
\end{lemma}

Then, continuous functionals of $\hat{S}_{1T}$ no longer converge to those of 
$S_{1}+\delta D_{1}$ under $H_{1T}$, but the estimation effect also has to
be taken into account using the following assumption. Let $Z\left( \Psi
\right) $ be a normal vector with zero mean and covariance matrix $\Psi .$

\begin{assumption}[Parameter estimation]

Under $H_{1T}$, the estimator $\hat{\theta}_{T}$ admits the asymptotic
linear expansion 
\begin{equation}
T^{1/2}\left( \hat{\theta}_{T}-\theta _{0}\right) =\delta \xi _{0}+\frac{1}{{%
T^{1/2}}}\sum_{t=1}^{T}\ell _{t}\left( Y_{t},\Omega _{t}\right) +o_{p}\left(
1\right) ,  \label{eq:estexpG}
\end{equation}%
where $\xi _{0}$ is a $m\times 1$ vector and the summands $\ell _{t}$
constitute a martingale difference sequence with respect to $\Omega _{t}$,
such that

\begin{description}
\item[(A)] $\mathop{\rm E}\nolimits \left[\ell_t \left(Y_t,\Omega
_{t}\right)\mid \Omega_t\right] =0$ and $T^{-1}\sum_{t=1}^{T} \mathop{\rm E}%
\nolimits \left[\ell_t \left(Y_t,\Omega _{t}\right) \ell_t \left(Y_t,\Omega
_{t}\right)^{\prime }\mid \Omega_t\right]\overset{p}{\rightarrow }\Psi . $

\item[(B)] Lindeberg condition 
$T^{-1}\sum_{t=1}^{T} \mathop{\rm E}\nolimits \left[\left \| \ell_t
\left(Y_t,\Omega _{t}\right)\right \|^2 \mathbbm{1}\! \left \{ {T^{-1/2}}
\left \| \ell_t \left(Y_t,\Omega _{t}\right)\right \|>\varepsilon \right \}
\mid \Omega _{t}\right]\overset{p}{\rightarrow }0 $ holds.

\item[(C)] There exists a finite function $W_{1}\left( u\right) $, such that 
$T^{-1}\sum_{t=1}^{T}\mathop{\rm E}\nolimits\left[ I_{t,\theta _{0}}\left(
u\right) \ell _{t}\left( Y_{t},\Omega _{t}\right) \mid \Omega _{t}\right]
\rightarrow _{p}W_{1}\left( u\right) $ uniformly in $u$.
\end{description}

In particular, under $H_{0}$, $\delta \xi _{0}=0$, the estimate $\hat{\theta}%
_{T}$ is centered and $T^{1/2}\left( \hat{\theta}_{T}-\theta _{0}\right) $
converges in distribution to $Z\left( \Psi \right) $.\newline
\end{assumption}

Assumption 4(A) and 4(B) hold for the MLE of many popular discrete
models, including dynamic probit and logit and general discrete choice
models. As an example consider estimates $\hat{\theta}_{T}$, which
are asymptotically equivalent to the (conditional) maximum likelihood
estimates, i.e.,
\begin{equation*}
T^{1/2}\left( \hat{\theta}_{T}-\theta _{0}\right) =-\frac{B_{0}^{-1}}{T^{1/2}%
}\sum_{t=1}^{T}s_{t}\left( Y_{t},\Omega _{t}\right) +o_{p}\left( 1\right) ,
\end{equation*}%
where $s_{t}\left( k,\Omega _{t}\right) :=\dot{f}_{{t,\theta _{0}}%
}\left( k\mid \Omega _{t}\right) /f_{{t,\theta _{0}}}\left( k\mid \Omega
_{t}\right) $ is the score function and $B_{0}$ is a
symmetric $m\times m$ positive definite matrix given by the limit of
the Hessian, 
\begin{equation*}
B_{0}:=\mathop{\rm plim}_{T\rightarrow \infty }\frac{1}{T}%
\sum_{t=1}^{T}\sum_{k=1}^{K}s_{t}\left( k,\Omega _{t}\right) \dot{f}_{{%
t,\theta _{0}}}\left( k\mid \Omega _{t}\right) ^{\prime }.
\end{equation*}%
Under $H_{1T}$, $\mathop{\rm E}\nolimits \left[
s_{t}\left( Y_{t},\Omega _{t}\right) \mid \Omega _{t}\right] =\delta
T^{-1/2}\sum_{k=1}^{K}s_{t}\left( k,\Omega _{t}\right) h_{t}\left( k\mid
\Omega _{t}\right) $. Then equation (\ref{eq:estexpG}) holds with $%
\xi _{0}=-\mathop{\rm plim}_{T\rightarrow \infty
}B_{0}^{-1}T^{-1}\sum_{t=1}^{T}\sum_{k=1}^{K}s_{t}\left( k,\Omega
_{t}\right) h_{{t}}\left( k\mid \Omega _{t}\right) $ and \newline
$\ell _{t}\left( Y_{t},\Omega _{t}\right) =-B_{0}^{-1}s_{t}\left(
Y_{t},\Omega _{t}\right) +B_{0}^{-1}\sum_{k=1}^{K}s_{t}\left( k,\Omega
_{t}\right) h_{{t}}\left( k\mid \Omega _{t}\right)$.\bigskip

We can derive the covariance matrix between the process $S_{1T}\left(
u\right) $ and $T^{1/2}\left( \hat{\theta}_{T}-\theta _{0}\right) $ and
obtain joint convergence results, so under $H_{1T}$%
\begin{equation}
\left( S_{1T},T^{1/2}\left( \hat{\theta}_{T}-\theta _{0}\right) \right)
\Rightarrow \left( S_{1\infty }+\delta D_{1},Z\left( \Psi \right) +\delta
\xi _{0}\right) ,  \label{eq:pujoint}
\end{equation}%
where the covariance function between $S_{1\infty }$ and $Z\left( \Psi
\right) $ is $W_{1}(u)$.

We can state now the result on the asymptotic distribution of the empirical
process~$\hat{S}_{1T}$ under local alternatives, whose drift is different
with respect to the case without estimated parameters.

\begin{theorem}
\label{thm:limitNormSLAPEE} Suppose Assumptions 1-4 hold. Under $H_{1T},$%
\begin{equation*}
\hat{S}_{1T}\Rightarrow \hat{S}_{1\infty} + \delta \left \{ D_1+\xi
_{0}^{\prime }L_1 \right \},  \label{eq:SLAPEE}
\end{equation*}%
where $\hat{S}_{1\infty}:=S_{1\infty}+Z\left( \Psi \right) ^{\prime }L_1 $
is a Gaussian process with zero mean and variance function $V_1\left(
u,v\right) +L_1\left( u\right)^{\prime } \Psi L_1\left( v\right) +W_1\left(
u\right) ^{\prime }L_1\left( v\right) +W_1\left( v\right) ^{\prime }
L_1\left( u\right) $.
\end{theorem}

\section{EMPIRICAL PROCESSES FOR DYNAMIC SPECIFICATION}

Test statistics based on ${S}_{1T}$, ${R}_{1T}$ and ${R}_{1T,M}$ verify that
the conditional distribution of $Y_{t}$ is right on \textit{average} across
all possible $\Omega _{t}$, so these tests might not capture all sources of
misspecification. This issue is raised in Corradi and Swanson (2006),
Delgado and Stute (2008) and Kheifets (2015) in relation to testing
continuous distributions. However, it is not possible to develop
specification tests conditioned on infinite dimensional values of $\Omega
_{t}$. Instead of truncating $\Omega _{t}$ or restricting the class of
models, we consider $S_{2T}$, a biparameter analog of $S_{1T}$ to control
the possible dynamic misspecification. From Lemma~\ref{lem:mds}, since under 
$H_{0}$, $I_{t,\theta _{0}}\left( u_{1}\right) -u_{1}$ is a MDS, $%
I_{t,\theta _{0}}\left( u_{1}\right) I_{t-1,\theta _{0}}\left( u_{2}\right)
-u_{1}u_{2}$ is centered around zero, and moreover%
\begin{equation*}
\mathop{\rm E}\nolimits\left[ I_{t,\theta _{0}}\left( u_{1}\right)
I_{t-1,\theta _{0}}\left( u_{2}\right) \mid \Omega _{t-1}\right]
=u_{1}u_{2},\ \ \ a.s.
\end{equation*}%
This motivates us to develop tests based on $S_{2T}$ defined in (\ref{eq:S2}). This process also has zero mean under the null and identifies not only
departures from the null derived from deviations of the unconditional
expectation of $I_{t,\theta _{0}}\left( u\right) $ from $u,$ but also from a
possible failure of the martingale property, so that $I_{t,\theta
_{0}}\left( u_{1}\right) $ and $I_{t-1,\theta _{0}}\left( u_{2}\right) $
would become correlated. This idea is similar to that exploited in Kheifets'
(2015) in the context of conditional distribution testing for continuous
distributions, where different methods of checking the independence property
of the PIT are proposed. Alternative statistics exploiting the lack of
correlations with any other lag could be proposed, but we expect that low
lags are typically more useful for detecting general forms of
misspecification.

One could also consider a biparameter analog of ${R}_{1T,M}$, i.e. for some $%
M=1,2,\ldots,$ 
\begin{equation*}
{R}_{2T,M}\left( u\right) :=\frac{1}{(T-1)^{1/2}M}\sum_{t=2}^{T}%
\sum_{m=1}^{M}\left(\mathbbm{1}\! \left \{ U_{t,m}^{r}\left( {\theta}%
_{0}\right) \leq u_1\right \} \mathbbm{1}\! \left \{ U_{t-1,m}^{r}\left( {%
\theta}_{0}\right) \leq u_2\right \} -u_1 u_2\right),
\end{equation*}
where $u=\left(u_1,u_2\right) \in [0,1]^2$. In particular, a bivariate
analog of ${R}_{1T}$, ${R}_{2T}\left( u\right) :={R}_{2T,1}\left( u\right)$,
is introduced in Kheifets and Velasco (2013). Tests based on ${R}_{2T}$ and $%
{R}_{2T,M}$ involve randomized transforms and therefore suffer from power
loss compared to tests based on the nonrandomized transform.

Note, that $S_{2T}\left( u\right) -u_{1}S_{1T-1}\left( u_{2}\right) $ is a
martingale. This observation will allow us to derive weak convergence of $%
S_{2T}$ by employing limiting theorems for MDS. Properties of $R_{2T}$ were
established in Kheifets and Velasco (2013) and could be extended to $%
R_{2T,M} $. Here we discuss the properties of $S_{2T}$ when we estimate $%
\theta _{0}.$

In practice we use the process 
\begin{equation*}
\hat{S}_{2T}\left( u\right) :=\frac{1}{(T-1)^{1/2}}\sum_{t=2}^{T}\left\{
I_{t,\hat{\theta}_{T}}\left( u_{1}\right) I_{t,\hat{\theta}_{T}}\left(
u_{2}\right) -u_{1}u_{2}\right\} ,
\end{equation*}%
where we can write under $H_{1T}$ 
\begin{equation}
\hat{S}_{2T}\left( u\right) =S_{2T}\left( u\right) +T^{1/2}\left( \hat{\theta%
}_{T}-\theta _{0}\right) ^{\prime }\frac{1}{T}\sum_{t=2}^{T}\nabla
_{2,t}\left( u\right) +o_{p}\left( 1\right) ,  \label{eq:puexpansion2}
\end{equation}%
uniformly in $u$, where \newline
$\nabla _{2,t}\left( u\right) :=I_{t-1,{\theta }_{0}}\left( u_{2}\right)
\nabla \left( F_{t,\theta _{0}}\left( \cdot \mid \Omega _{t}\right)
,u_{1}\right) +u_{1}\nabla \left( F_{t-1,\theta _{0}}\left( \cdot \mid
\Omega _{t-1}\right) ,u_{2}\right) $ and the asymptotic covariance function
is $W_{2}\left( u\right) :=\mathop{\rm ACov}\nolimits\left( S_{2T}\left(
u\right) ,T^{1/2}\left( \hat{\theta}_{T}-\theta _{0}\right) \right) $. To
study the asymptotic properties of the biparameter process we introduce the
next assumption, which extends Assumption~2.

\begin{assumption}
Under $H_{1T}$, there exist finite functions $D_{2}\left( u\right) $ and $%
L_{2}\left( u\right) $, such that uniformly in~$u$

\begin{description}
\item[(A)] 
$T^{-1}\sum_{t=2}^{T}\left\{ I_{t-1,\theta _{0}}\left( u_{2}\right) d\left(
H_{t}\left( \cdot \mid \Omega _{t}\right) ,F_{t,\theta _{0}}\left( \cdot
\mid \Omega _{t}\right) ,u_{1}\right) \right. $ \newline
$+\left. u_{1}d\left( H_{t}\left( \cdot \mid \Omega _{t}\right) ,F_{t,\theta
_{0}}\left( \cdot \mid \Omega _{t}\right) ,u_{2}\right) \right\} \rightarrow
_{p}D_{2}\left( u\right) $. 

\item[(B)] 
$T^{-1}\sum_{t=2}^{T}\nabla_{2,t}\left(u\right)\to_p L_2\left( u\right) $. 
\end{description}
\end{assumption}

Note that the second terms in the definitions of $D_{2}$ and $L_{2}$
correspond to $u_{1}D_{1}(u_{2})$ and $u_{1}L_{1}(u_{2})$ respectively, the
equivalent for the single parameter process $S_{1T}$, but the first ones are
new. To state the next result, we need to assume existence of probabilistic
limits of several random functions. For the sake of presentation, we defer
precise statements to the Appendix, see Assumption A.

\begin{theorem}
\label{thm:limitNormSLAPEE2} Suppose that in addition to the conditions of
Theorem~\ref{thm:limitNormSLAPEE}, Assumption~5 and Assumption~A from the
Appendix hold. Under $H_{1T}$,%
\begin{equation*}
S_{2T}\left( u\right) \Rightarrow S_{2\infty}+\delta D_2,
\end{equation*}%
where $S_{2\infty}$ is a Gaussian process in $\left[ 0,1\right] $ with mean
zero and covariance function $V_{2\infty}\left( u,v\right)$ defined in the
Appendix. Under $H_{1T}$, if parameters are estimated,%
\begin{equation*}
\hat{S}_{2T}\Rightarrow \hat{S}_{2\infty}+\delta \left \{ D_2+\xi
_{0}^{\prime }L_2\right \},  \label{eq:SLAPEE2}
\end{equation*}%
where $\hat{S}_{2\infty}:=S_{2\infty}+Z\left(\Psi \right) ^{\prime }L_2$ is
a Gaussian process with zero mean and variance function $V_{2\infty}\left(
u,v\right) +L_2\left( u\right)^{\prime } \Psi L_2\left( v\right) +W_2\left(
u\right) ^{\prime }L_2\left( v\right) +W_2\left( v\right) ^{\prime }
L_2\left( u\right)$.
\end{theorem}

When $G_{t}\left( \cdot \mid \Omega _{t}\right) $ is different from $%
F_{t,\theta _{0}}\left( \cdot \mid \Omega _{t}\right) $ such that $D_{2}$ is
non-zero, the test based on $\hat{S}_{2T}$ has nontrivial power in the
direction of $H_{1T}$. In contrast to the univariate case with~$S_{1T}$, the
first term in the definition of $D_{2}$ contains correlation with the past
information and can therefore capture dynamic misspecification when 
this induces in such a correlation, even if the unconditional
expectation of $d$, which appears in the second term $u_{1}D_{1}(u_{2})$, 
is zero. This fact is crucial if misspecification occurs in the dynamics and
not only in the link function or other static aspects of the model.

\section{BOOTSTRAP TESTS}

To test $H_{0}$ we consider Cramer-von Mises, Kolmogorov-Smirnov or any
other continuous functionals of $\hat S_{jT}$, $j=1,2$, $\eta \left( \hat
S_{jT}\right)$. Then consistency properties of specification tests based on $%
\hat S_{jT}$ can be derived using the discussion in the previous sections by
applying the continuous mapping theorem, so we omit the proof of the
following result.

\begin{theorem}
\label{thm:limitNormSLA} Suppose that conditions of Theorem~\ref%
{thm:limitNormSLAPEE2} hold. Under $H_{1T}$,%
\begin{equation*}
\eta \left( \hat{S}_{jT}\right) \rightarrow _{d}\eta \left( \hat{S}_{j\infty
}\right) ,\ \ j=1,2.
\end{equation*}
\end{theorem}

Since the asymptotic distributions of $S_{jT}\left( u\right)$ are model
dependent, and those of $\hat{S}_{jT}\left( u\right)$ further depend on the
estimation effect, we need to resort to bootstrap methods to implement our
tests in practice. In the literature, there are several resampling methods
suitable for dependent data, but since under $H_{0}$ the parametric
conditional distribution is fully specified, we apply a conditional
parametric bootstrap algorithm that only requires to make draws from $%
F_{t,\hat \theta}\left(\cdot \mid \Omega _{t}\right)$ to mimic the null
distribution of the test statistics. For a discussion of the parametric
bootstrap see Stute et al.\ (1993) and Andrews (1997), which can be adapted to the complications with
information truncation and initialization arising in the dynamic case using
the discussion in Bai (2003).

To estimate the true $1-\alpha $ quantiles $c_j\left( \theta _{0}\right) $
of the null asymptotic distribution of the test statistics, given by some
continuous functional $\eta$ applied to $\hat{S}_{j\infty}$ with $\delta =0$%
, we implement the following steps.

\begin{enumerate}
\item Estimate the model with data $\left( Y_{t},X^{\prime }_{t}\right) $, $%
t=1,2,...,T$, get parameter estimator $\hat{\theta}_{T}$ and compute test
statistics $\eta (\hat{S}_{jT})$.

\item Simulate $Y_{t}^{\ast }$ with $F_{\hat{\theta_T}}\left(\cdot \mid \Omega
_{t}^{\ast }\right)$ recursively for $t=1,2,...,T$, where the bootstrap
information set is $\Omega _{t}^{\ast }=\left(X_{t}, Y_{t-1}^{\ast
},X_{t-1},Y_{t-2}^{\ast },X_{t-2},...\right) $.

\item Estimate the model with simulated data $Y_{t}^{\ast }$, get $\hat{%
\theta}_{T}^{\ast }$ using the same method as for $\hat{\theta}_{T},$ get
bootstrapped test statistics $\eta \left( \hat{S}_{jT}^{\ast }\right) $.

\item Repeat 2-3 $B$ times, compute the percentiles of the empirical
distribution of the $B$ bootstrapped test statistics.

\item Reject $H_{0}$ if $\eta \left( \hat{S}_{jT}\right)$ is greater than
the $(1-\alpha )$th percentile of the empirical distribution of the $B$
bootstrapped test statistics denoted by $\hat{c}^{\ast}_{jB}\left( \hat{%
\theta}_{T}\right) $.
\end{enumerate}

To analyze the properties of our parametric bootstrap, we need to assume
that the same conditions on the estimation method hold for both for original
and resampled data. More formally, we have

\begin{assumption}
\begin{description}
\item[(A)] The conditional distribution of $Y_t$ conditional on $\Omega_t$
coincides with the conditional distribution of $Y_t$ conditional on $%
\Omega_t\cup \{X^{\prime }_k\}_{k=t+1}^T$.

\item[(B)] Suppose that the sample is generated by $F_{\theta_T}$, for some
nonrandom sequence $\theta_T$ converging to $\theta_0$, i.e.\ we have a
triangular array of random variables $\{Y_{Tt}:t=1,2,\ldots,T\}$ with $(T,t)$
element generated by $F_{\theta _{T}}(\cdot \mid \Omega _{Tt})$, where 
\newline 
$\Omega _{Tt}=\left \{X_{t},Y_{Tt-1},X_{t-1},Y_{Tt-2},X_{t-2},\ldots \right \}
$. Then the estimator $\hat \theta_T$ of $\theta_T$ admits an asymptotic
linear expansion as in Assumption~4. Moreover, assume that under the
alternative $H_{1},$ there exists some $\theta _{1}\in\Theta$ so that $\theta _{1}=%
\mathop{\rm
plim}_{T\rightarrow \infty }\hat{\theta}_{T}.$
\end{description}
\end{assumption}

This assumption insures that by simulating from the conditional
distribution $F_{\theta_T}$ we obtain the correct joint distribution
of $S_{jT}$ and $T^{1/2}\left(\hat\theta_T-\theta_T\right)$ in
parallel to those required in Theorems 1-2. Assumption 6 (A) says
that $Y_{t}$ and future $X_{t}$ are independent
conditionally on past information, i.e.\ that there is no direct feedback
effect. For example, in a latent variable form of the ordered probit model,
this assumption translates to strict exogeneity, i.e. that innovations are
independent of future $X_{t}$. Dependence between $Y_{t}$ and future
$X_{t}$ is still allowed through serial dependence in $X_{t}$ and
$Y_{t}$. Assumption~6 (B) is similar to Condition (5.5) in Burke et
al.\ (1979), Assumption (A1) in Stute et
al.\ (1993) and Assumption~E2 in Andrews (1997), and introduces a
triangular array version of the expansion and central limit theorem for parameter estimates, see also the discussion in Section~4.1 in Andrews (1997).

We obtain the following result.

\begin{theorem}
\label{thm:limitNormSBoot} Suppose that in addition to conditions of Theorem~%
\ref{thm:limitNormSLAPEE2}, Assumption~6 holds. Under $H_{1T},$ as $B,T$ $%
\rightarrow \infty ,$%
\begin{equation*}
\eta \left( \hat{S}_{jT}^{\ast }\right) \rightarrow _{d}\eta \left( \hat{S}%
_{j\infty }\right) ,\qquad j=1,2,
\end{equation*}%
in probability, so $\hat{c}_{jB}^{\ast }\left( \hat{\theta}_{T}\right)
\rightarrow _{p}c_{j}\left( \theta _{0}\right) $, and therefore, under $%
H_{0},$ $\Pr \left( \eta \left( \hat{S}_{jT}\right) >\hat{c}_{jB}^{\ast
}\left( \hat{\theta}_{T}\right) \right) \rightarrow \alpha $. Suppose also that the conditions of Theorem~%
\ref{thm:limitNormSLAPEE2} hold for any $\theta_0 \in\Theta$. Under $H_{1},$
as $B,T$ $\rightarrow \infty ,\ \hat{c}_{jB}^{\ast }\left( \hat{\theta}%
_{T}\right) =O_{p}\left( 1\right) $.
\end{theorem}

This theorem shows that the bootstrap test statistic has the same limit
distribution as the original one under local alternatives, so that under the
null we get the right asymptotic size using bootstrap estimated critical
values and that under local alternatives we get non trivial power when the
drifts of the stochastic processes $\hat{S}_{1T}$ and $\hat{S}_{2T}$ are non
negligible. Similarly, under fixed alternatives we are able to get a
bootstrap consistent test when the asymptotic test is consistent itself,
i.e. $\lim_{T\rightarrow \infty }\Pr \left( \eta \left( \hat{S}_{jT}\right) >%
\hat{c}_{jB}^{\ast }\left( \hat{\theta}_{T}\right) \right) =1$ if $\eta
\left( \hat{S}_{jT}\right) $ diverges asymptotically.

\section{APPLICATION AND SIMULATIONS}

In this section we use a Monte Carlo simulation exercise to investigate the
finite sample properties of the tests proposed in this paper. We take as
reference the dynamic ordered discrete choice models investigated in Basu
and de Jong (2007) for the modeling of the monetary policy conducted by the
Federal Reserve (FED). The dependent variable uses the following
codification of the changes in the reference interest rate in US, the
federal funds rate $i_{t}$,%
\begin{equation*}
Y_{t}=\left\{ 
\begin{array}{ccc}
1 & \ \text{if} & \Delta i_{t}<-0.25 \\ 
2 & \ \text{if} & -0.25\leq \Delta i_{t}<0 \\ 
3 & \ \text{if} & 0\leq \Delta i_{t}<0.25 \\ 
4 & \ \text{if} & \Delta i_{t}\geq 0.25.%
\end{array}%
\right. 
\end{equation*}

Data is monthly and spans January 1990 to December 2006, leading to $T=204$
complete observations. The explanatory variables that Basu and de Jong
(2007) used to explain the decisions of the FED on $\Delta i_{t}$ are the
current value and 4 lags of inflation $\left( \inf \right) $, the current
value and a lag of four different measures of output gap $\left( out\right) $
and a series of dummies that describe the decision of the FED in the
previous period, $dum1_t=I(\Delta i_{t-1}<0),\ dum2_t=I(\Delta i_{t-1}>0),\
dum3_t=I(\Delta i_{t-1}<-0.25),\ dum4_t=I(\Delta i_{t-1}>0.25).$ Instead of
these four dummies, we implement an AR$\left( 1\right) $, 'dynamic' version
with one lag of the discrete $Y_{t}$ as explanatory variable (and a version
without lags that we refer to as 'static' to serve as a benchmark to the
inclusion of lagged endogenous variables in $\Omega _{t})$. We consider both
the Logit and Probit versions of the models. We fit four versions of the
basic model based on different definitions of the output gap and conditional
on the series of inflation and output gap and on the parameter estimates
obtained, we simulate series $Y_{t}$ and conduct our tests on these (see
Monte Carlo scenarios in Table~\ref{t:rejscenarios}). 

\begin{table}
\begin{center}
\caption{Scenarios for Monte Carlo simulations.}
\label{t:rejscenarios}
\begin{tabular}{cl}
\hline
\hline
 Scenario & Null and Alternative \tabularnewline
\hline
Size 1  &  $H_0:$ static probit\tabularnewline
Size 2  &  $H_0:$ static logit\tabularnewline
Power 1  &  $H_0:$ static probit vs  $H_1:$ static logit\tabularnewline
Power 2  &  $H_0:$ static probit vs  $H_1:$ dynamic probit\tabularnewline
Power 3  &  $H_0:$ static probit vs  $H_1:$ dynamic logit\tabularnewline
\hline
\end{tabular}
\end{center}
\end{table}

The four choices of output gap lead to Models I-IV. The output gap is the
percentage deviation of the actual from the potential output, which is
interpolated to obtain a series of monthly frequency by replicating the GDP
observation for any quarter to all the months in that quarter. Then two
different measures of potential output are used: the potential output series
provided by the Congressional Budget Office and a potential output series
constructed in a real-time setting using the HP filter, leading to Models I
and II. Apart from output gap, other measures of economic activity are used,
such as unemployment rate and capacity utilization, leading to Models III
and IV. Data sources are described in Basu and de Jong (2007).

We compare the performance of our tests with an alternative test which is
also omnibus and does not require smoothing (and choice of smoothing
parameters). Two general approaches can be adapted to our setup: the test of
the Generalized Linear Model (GLM) of Stute and Zhu (2002) and the
Conditional Kolmogorov test of Andrews (1997), as discussed in Mora and
Moro-Egido (2007). The first one is a test based on a marked empirical
process for testing the null $H_{0}^{\prime }:\quad \mathop{\rm E}\nolimits%
\left[ Y\mid \tilde{X}=x\right] =m_{\tilde{\beta}_{01}}\left( x^{^{\prime }}%
\tilde{\beta}_{02}\right) $, where $m_{\tilde{\beta}_{01}}(\cdot )$ is a
parametric link function and $\tilde{\beta}_{01},\tilde{\beta}_{02}$ are
finite dimensional parameters. In the cases where $Y$ takes only two values $%
\{0,1\}$, the conditional mean coincides with the conditional probability of 
$Y=1$ and the null is similar to our $H_{0}$ if we were considering an i.i.d
setup. To test $Y_{t}\mid \tilde{X}_{t}\sim P_{\tilde{\beta}_{01}}\left(
\cdot \mid \tilde{X}_{t}^{^{\prime }}\tilde{\beta}_{20}\right) $ define the
process 
\begin{equation*}
{Z}_{T}\left( y\right) :=\frac{1}{T^{1/2}}\sum_{t=1}^{T}\mathbbm{1}\!\left\{ 
\tilde{X}_{t}^{^{\prime }}\tilde{\beta}_{20}\leq y\right\} \left[ Y_{t}-P_{%
\tilde{\beta}_{01}}\left( Y_{t}=1\mid \tilde{X}_{t}^{^{\prime }}\tilde{\beta}%
_{20}\right) \right] ,\qquad y\in \mathop{\mathbb R}\nolimits.
\end{equation*}%
The second test by Andrews is obtained by substituting $\mathbbm{1}\!\left\{ 
\tilde{X}_{t}^{^{\prime }}\tilde{\beta}\leq y\right\} $ with $\mathbbm{1}%
\!\left\{ \tilde{X}_{t}\leq \tilde{x}\right\} $ (where $\tilde{x}$ is a real
vector of dimension of $\tilde{X}_{t}$) in ${Z}_{T}$, but since it always
underperforms according to simulations of Mora and Moro-Egido (2007), it is
not considered here. If $Y$ takes values $\{1,\ldots ,K\}$, Mora and
Moro-Egido (2007) substitute testing $H_{0}$ by $K$ tests of the hypotheses $%
Y_{jt}\mid \tilde{X}_{t}\sim P_{j,\tilde{\beta}_{01}}\left( Y_{t}\mid \tilde{%
X}_{t}^{^{\prime }}\tilde{\beta}_{20}\right) $, with corresponding processes 
${Z}_{j,T}$, where $Y_{jt}=\mathbbm{1}\!\left\{ Y_{t}=j\right\} $ and $%
j=1,2,\ldots ,K$, then the resulting pooled test statistics are 
\begin{equation*}
\eta _{Z}^{CvM}=T^{-1}\sum_{j=1}^{K}\sum_{t=1}^{T}{Z}_{j,T}\left( \tilde{X}%
_{t}^{^{\prime }}\tilde{\beta}_{20}\right) ^{2}
\end{equation*}%
and 
\begin{equation*}
\eta _{Z}^{KS}=T^{-1}\max_{j=1,\ldots ,K}\sum_{t=1}^{T}{Z}_{j,T}\left( 
\tilde{X}_{t}^{^{\prime }}\tilde{\beta}_{20}\right) ^{2},
\end{equation*}%
which we call the CvM and KS tests respectively. To apply these tests to our
model, let $\tilde{X}_{t}=\left( X_{t}^{^{\prime }},Y_{t-1},1\right) ^{\prime }
$ and $\tilde{\beta}=\left( \beta ^{^{\prime }},\rho,-\tau_1 \right) ^{\prime }$
and take the corresponding link functions.

We analyze tests based on $S_{1T}$, $R_{1T,M}$, $R_{1T}$ and $S_{2T}$, $%
R_{2T,M}$, $R_{2T}$ and $Z_T$. In all cases we use Kolmogorov-Smirnov (KS)
and Cramer-von Mises (CvM) measures. We only consider feasible bootstrap
versions of tests based on $\hat S_{1T}$, $\hat R_{1T,M}$, etc, where we
replace $\theta _{0}$ by root-$T\ $consistent estimates $\hat{\theta}_{T}$,
the ML estimator in our case. We are not aware of any theoretical results
for bootstrap assisted tests based on $\hat Z_T$ in our setup, although Mora
and Moro-Egido (2007) provide some simulations.

\renewcommand*{\arraystretch}{.75}
\begin{table}\begin{center}
\caption{ML estimates and standard errors of Models I-IV with static and dynamic specifications and Probit link function applied to the real US data, $T=204$.}
\label{t:estprobit}
\setlength{\tabcolsep}{1pt}
\begin{tabular}{ccccccccccc}\hline\hline
 & I-static  & I-dynamic  & II-static  & II-dynamic  & III-static  & III-dynamic  & IV-static  & IV-dynamic \tabularnewline\hline
  $\tau_1$  &  $ -4.81 $   &  $ -2.07 $   &  $ -3.31 $   &  $ -1.05 $   &  $ -3.15 $   &  $ -1.17 $   &  $ -3.41 $   &  $ -1.48 $  \tabularnewline
 &  $ (0.51) $  &  $ (0.66) $  &  $ (0.35) $  &  $ (0.47) $  &  $ (0.36) $  &  $ (0.48) $  &  $ (0.37) $  &  $ (0.50) $ \tabularnewline
  $\tau_{2}$  &  $ -4.05 $   &  $ -1.14 $   &  $ -2.64 $   &  $ -0.19 $   &  $ -2.34 $   &  $ -0.20 $   &  $ -2.57 $   &  $ -0.50 $  \tabularnewline
 &  $ (0.47) $  &  $ (0.64) $  &  $ (0.31) $  &  $ (0.46) $  &  $ (0.32) $  &  $ (0.47) $  &  $ (0.32) $  &  $ (0.48) $ \tabularnewline
  $\tau_{3}$  &  $ -1.72 $   &  $ 1.66 $   &  $ -0.39 $   &  $ 2.60 $   &  $ 0.09 $   &  $ 2.62 $   &  $ -0.11 $   &  $ 2.29 $  \tabularnewline
 &  $ (0.40) $  &  $ (0.63) $  &  $ (0.26) $  &  $ (0.48) $  &  $ (0.28) $  &  $ (0.48) $  &  $ (0.27) $  &  $ (0.49) $ \tabularnewline
  $inf$  &  $ -1.39 $   &  $ -1.36 $   &  $ -1.51 $   &  $ -1.60 $   &  $ -1.83 $   &  $ -1.82 $   &  $ -1.70 $   &  $ -1.70 $  \tabularnewline
 &  $ (0.68) $  &  $ (0.72) $  &  $ (0.67) $  &  $ (0.71) $  &  $ (0.69) $  &  $ (0.73) $  &  $ (0.69) $  &  $ (0.73) $ \tabularnewline
  $inf_{-1}$  &  $ 1.86 $   &  $ 2.90 $   &  $ 1.94 $   &  $ 3.05 $   &  $ 2.05 $   &  $ 3.07 $   &  $ 2.14 $   &  $ 3.01 $  \tabularnewline
 &  $ (0.99) $  &  $ (1.06) $  &  $ (0.98) $  &  $ (1.06) $  &  $ (1.00) $  &  $ (1.07) $  &  $ (1.01) $  &  $ (1.07) $ \tabularnewline
  $inf_{-2}$  &  $ -1.30 $   &  $ -2.81 $   &  $ -1.27 $   &  $ -2.80 $   &  $ -1.60 $   &  $ -2.92 $   &  $ -2.12 $   &  $ -3.11 $  \tabularnewline
 &  $ (0.98) $  &  $ (1.07) $  &  $ (0.97) $  &  $ (1.06) $  &  $ (0.99) $  &  $ (1.07) $  &  $ (1.02) $  &  $ (1.09) $ \tabularnewline
  $inf_{-3}$  &  $ 1.39 $   &  $ 2.44 $   &  $ 1.60 $   &  $ 2.74 $   &  $ 1.79 $   &  $ 2.79 $   &  $ 1.27 $   &  $ 2.33 $  \tabularnewline
 &  $ (0.99) $  &  $ (1.06) $  &  $ (0.98) $  &  $ (1.06) $  &  $ (1.00) $  &  $ (1.08) $  &  $ (1.03) $  &  $ (1.09) $ \tabularnewline
  $inf_{-4}$  &  $ 0.43 $   &  $ -0.53 $   &  $ -0.23 $   &  $ -1.05 $   &  $ -0.00 $   &  $ -0.85 $   &  $ 0.88 $   &  $ -0.20 $  \tabularnewline
 &  $ (0.68) $  &  $ (0.73) $  &  $ (0.66) $  &  $ (0.71) $  &  $ (0.67) $  &  $ (0.73) $  &  $ (0.71) $  &  $ (0.76) $ \tabularnewline
  $out$  &  $ -1.02 $   &  $ -1.02 $   &  $ 0.36 $   &  $ 0.40 $   &  $ 3.35 $   &  $ 2.54 $   &  $ -0.98 $   &  $ -0.62 $  \tabularnewline
 &  $ (0.30) $  &  $ (0.33) $  &  $ (0.59) $  &  $ (0.63) $  &  $ (0.68) $  &  $ (0.74) $  &  $ (0.22) $  &  $ (0.23) $ \tabularnewline
  $out_{-1}$  &  $ 0.81 $   &  $ 0.90 $   &  $ 0.84 $   &  $ 0.65 $   &  $ 2.48 $   &  $ 0.95 $   &  $ -1.03 $   &  $ -0.65 $  \tabularnewline
 &  $ (0.29) $  &  $ (0.32) $  &  $ (0.59) $  &  $ (0.64) $  &  $ (0.67) $  &  $ (0.73) $  &  $ (0.22) $  &  $ (0.23) $ \tabularnewline
 $Y_{-1}$ &  ---   &  $ -1.08 $   &  ---   &  $ -1.12 $   &  ---   &  $ -1.03 $   &  ---   &  $ -0.94 $  \tabularnewline
 &       &  $ (0.15) $  &       &  $ (0.15) $  &       &  $ (0.16) $  &       &  $ (0.16) $ \tabularnewline\hline
\end{tabular}\end{center}
\end{table}

\begin{table}\begin{center}
\caption{ML estimates and standard errors of Models I-IV with static and dynamic specifications and Logit link function applied to the real US data, $T=204$.}
\label{t:estlogit}
\setlength{\tabcolsep}{1pt}
\begin{tabular}{ccccccccccc}\hline\hline
 & I-static  & I-dynamic  & II-static  & II-dynamic  & III-static  & III-dynamic  & IV-static  & IV-dynamic \tabularnewline\hline
  $\tau_1$  &  $ -8.46 $   &  $ -3.77 $   &  $ -6.01 $   &  $ -2.12 $   &  $ -5.61 $   &  $ -2.15 $   &  $ -6.15 $   &  $ -2.82 $  \tabularnewline
 &  $ (0.98) $  &  $ (1.20) $  &  $ (0.68) $  &  $ (0.83) $  &  $ (0.69) $  &  $ (0.85) $  &  $ (0.72) $  &  $ (0.89) $ \tabularnewline
  $\tau_{2}$  &  $ -7.03 $   &  $ -1.96 $   &  $ -4.71 $   &  $ -0.46 $   &  $ -4.12 $   &  $ -0.31 $   &  $ -4.56 $   &  $ -0.90 $  \tabularnewline
 &  $ (0.90) $  &  $ (1.17) $  &  $ (0.60) $  &  $ (0.81) $  &  $ (0.59) $  &  $ (0.83) $  &  $ (0.61) $  &  $ (0.86) $ \tabularnewline
  $\tau_{3}$  &  $ -3.00 $   &  $ 3.02 $   &  $ -0.85 $   &  $ 4.52 $   &  $ 0.07 $   &  $ 4.60 $   &  $ -0.24 $   &  $ 4.04 $  \tabularnewline
 &  $ (0.72) $  &  $ (1.12) $  &  $ (0.47) $  &  $ (0.84) $  &  $ (0.49) $  &  $ (0.86) $  &  $ (0.49) $  &  $ (0.87) $ \tabularnewline
  $inf$  &  $ -2.44 $   &  $ -2.29 $   &  $ -2.53 $   &  $ -2.89 $   &  $ -3.17 $   &  $ -3.28 $   &  $ -2.81 $   &  $ -3.06 $  \tabularnewline
 &  $ (1.21) $  &  $ (1.30) $  &  $ (1.21) $  &  $ (1.29) $  &  $ (1.21) $  &  $ (1.32) $  &  $ (1.22) $  &  $ (1.32) $ \tabularnewline
  $inf_{-1}$  &  $ 3.28 $   &  $ 4.95 $   &  $ 3.22 $   &  $ 5.46 $   &  $ 3.59 $   &  $ 5.43 $   &  $ 3.41 $   &  $ 5.31 $  \tabularnewline
 &  $ (1.78) $  &  $ (1.92) $  &  $ (1.77) $  &  $ (1.92) $  &  $ (1.76) $  &  $ (1.93) $  &  $ (1.82) $  &  $ (1.95) $ \tabularnewline
  $inf_{-2}$  &  $ -2.48 $   &  $ -5.02 $   &  $ -2.17 $   &  $ -5.22 $   &  $ -2.97 $   &  $ -5.21 $   &  $ -3.52 $   &  $ -5.40 $  \tabularnewline
 &  $ (1.74) $  &  $ (1.95) $  &  $ (1.73) $  &  $ (1.94) $  &  $ (1.76) $  &  $ (1.95) $  &  $ (1.86) $  &  $ (1.99) $ \tabularnewline
  $inf_{-3}$  &  $ 2.42 $   &  $ 4.36 $   &  $ 2.61 $   &  $ 5.20 $   &  $ 2.94 $   &  $ 5.11 $   &  $ 1.65 $   &  $ 4.02 $  \tabularnewline
 &  $ (1.75) $  &  $ (1.92) $  &  $ (1.75) $  &  $ (1.93) $  &  $ (1.77) $  &  $ (1.95) $  &  $ (1.86) $  &  $ (1.99) $ \tabularnewline
  $inf_{-4}$  &  $ 0.93 $   &  $ -0.87 $   &  $ -0.17 $   &  $ -1.88 $   &  $ 0.32 $   &  $ -1.54 $   &  $ 2.11 $   &  $ -0.28 $  \tabularnewline
 &  $ (1.20) $  &  $ (1.32) $  &  $ (1.18) $  &  $ (1.28) $  &  $ (1.19) $  &  $ (1.30) $  &  $ (1.27) $  &  $ (1.36) $ \tabularnewline
  $out$  &  $ -1.78 $   &  $ -1.79 $   &  $ 0.43 $   &  $ 0.63 $   &  $ 5.87 $   &  $ 4.12 $   &  $ -1.83 $   &  $ -1.15 $  \tabularnewline
 &  $ (0.54) $  &  $ (0.60) $  &  $ (1.04) $  &  $ (1.14) $  &  $ (1.24) $  &  $ (1.34) $  &  $ (0.40) $  &  $ (0.42) $ \tabularnewline
  $out_{-1}$  &  $ 1.43 $   &  $ 1.59 $   &  $ 1.61 $   &  $ 1.29 $   &  $ 4.21 $   &  $ 1.50 $   &  $ -1.88 $   &  $ -1.14 $  \tabularnewline
 &  $ (0.52) $  &  $ (0.59) $  &  $ (1.04) $  &  $ (1.15) $  &  $ (1.20) $  &  $ (1.33) $  &  $ (0.40) $  &  $ (0.42) $ \tabularnewline
 $Y_{-1}$ &  ---   &  $ -1.98 $   &  ---   &  $ -2.04 $   &  ---   &  $ -1.86 $   &  ---   &  $ -1.71 $  \tabularnewline
 &       &  $ (0.28) $  &       &  $ (0.27) $  &       &  $ (0.28) $  &       &  $ (0.28) $ \tabularnewline\hline
\end{tabular}\end{center}
\end{table}

\renewcommand*{\arraystretch}{1}

\begin{table}\begin{center}
\caption{P-values of Cramer -- von Misses tests for static Probit and Logit link function applied to the real US data, $T=204$.}
\label{t:pvalcvm}
\footnotesize\begin{tabular}{cccccccccccc}
\hline
\hline
 &  &  & $\hat S_{2T} $  & $\hat R_{2T,50} $  & $\hat R_{2T,25} $  & $\hat R_{2T} $  & $\hat S_{1T} $  & $\hat R_{1T,50} $  & $\hat R_{1T,25} $  & $\hat R_{1T} $  & $\hat Z_{T} $ \tabularnewline
\hline\multicolumn{11}{c}{$H_0:$ static probit} \tabularnewline
 & Model I &  & $ 0.001 $  & $ 0.001 $  & $ 0.001 $  & $ 0.237 $  & $ 0.009 $  & $ 0.026 $  & $ 0.078 $  & $ 0.516 $  & $ 0.244 $
\tabularnewline
 & Model II &  & $ 0.001 $  & $ 0.001 $  & $ 0.001 $  & $ 0.166 $  & $ 0.077 $  & $ 0.057 $  & $ 0.229 $  & $ 0.167 $  & $ 0.022 $
\tabularnewline
 & Model III &  & $ 0.001 $  & $ 0.001 $  & $ 0.001 $  & $ 0.307 $  & $ 0.492 $  & $ 0.632 $  & $ 0.616 $  & $ 0.731 $  & $ 0.109 $
\tabularnewline
 & Model IV &  & $ 0.001 $  & $ 0.002 $  & $ 0.002 $  & $ 0.496 $  & $ 0.721 $  & $ 0.509 $  & $ 0.582 $  & $ 0.668 $  & $ 0.268 $
\tabularnewline
\hline\multicolumn{11}{c}{$H_0:$ static logit} \tabularnewline
 & Model I &  & $ 0.001 $  & $ 0.001 $  & $ 0.001 $  & $ 0.152 $  & $ 0.021 $  & $ 0.079 $  & $ 0.221 $  & $ 0.793 $  & $ 0.199 $
\tabularnewline
 & Model II &  & $ 0.001 $  & $ 0.001 $  & $ 0.001 $  & $ 0.112 $  & $ 0.113 $  & $ 0.155 $  & $ 0.459 $  & $ 0.240 $  & $ 0.032 $
\tabularnewline
 & Model III &  & $ 0.001 $  & $ 0.001 $  & $ 0.001 $  & $ 0.360 $  & $ 0.314 $  & $ 0.493 $  & $ 0.541 $  & $ 0.745 $  & $ 0.171 $
\tabularnewline
 & Model IV &  & $ 0.001 $  & $ 0.001 $  & $ 0.001 $  & $ 0.448 $  & $ 0.890 $  & $ 0.804 $  & $ 0.899 $  & $ 0.634 $  & $ 0.272 $ \tabularnewline
\hline
\end{tabular}\end{center}
\end{table}

\begin{table}\begin{center}
\caption{P-values of Kolmogorov -- Smirnov tests for static Probit and Logit link function applied to the real US data, $T=204$.}
\label{t:pvalks}
\footnotesize\begin{tabular}{cccccccccccc}
\hline
\hline
 &  &  & $\hat S_{2T} $  & $\hat R_{2T,50} $  & $\hat R_{2T,25} $  & $\hat R_{2T} $  & $\hat S_{1T} $  & $\hat R_{1T,50} $  & $\hat R_{1T,25} $  & $\hat R_{1T} $  & $\hat Z_{T} $ \tabularnewline
\hline\multicolumn{11}{c}{$H_0:$ static probit} \tabularnewline
 & Model I &  & $ 0.003 $  & $ 0.002 $  & $ 0.002 $  & $ 0.082 $  & $ 0.047 $  & $ 0.193 $  & $ 0.372 $  & $ 0.354 $  & $ 0.392 $
\tabularnewline
 & Model II &  & $ 0.001 $  & $ 0.001 $  & $ 0.002 $  & $ 0.586 $  & $ 0.351 $  & $ 0.426 $  & $ 0.626 $  & $ 0.450 $  & $ 0.107 $
\tabularnewline
 & Model III &  & $ 0.001 $  & $ 0.001 $  & $ 0.001 $  & $ 0.155 $  & $ 0.454 $  & $ 0.435 $  & $ 0.244 $  & $ 0.742 $  & $ 0.124 $
\tabularnewline
 & Model IV &  & $ 0.001 $  & $ 0.002 $  & $ 0.002 $  & $ 0.799 $  & $ 0.936 $  & $ 0.913 $  & $ 0.801 $  & $ 0.355 $  & $ 0.230 $
\tabularnewline
\hline\multicolumn{11}{c}{$H_0:$ static logit} \tabularnewline
 & Model I &  & $ 0.001 $  & $ 0.001 $  & $ 0.001 $  & $ 0.133 $  & $ 0.010 $  & $ 0.050 $  & $ 0.212 $  & $ 0.684 $  & $ 0.220 $
\tabularnewline
 & Model II &  & $ 0.001 $  & $ 0.001 $  & $ 0.001 $  & $ 0.354 $  & $ 0.114 $  & $ 0.201 $  & $ 0.319 $  & $ 0.416 $  & $ 0.058 $
\tabularnewline
 & Model III &  & $ 0.001 $  & $ 0.001 $  & $ 0.001 $  & $ 0.149 $  & $ 0.511 $  & $ 0.472 $  & $ 0.350 $  & $ 0.642 $  & $ 0.173 $
\tabularnewline
 & Model IV &  & $ 0.002 $  & $ 0.002 $  & $ 0.001 $  & $ 0.769 $  & $ 0.975 $  & $ 0.968 $  & $ 0.867 $  & $ 0.411 $  & $ 0.207 $ \tabularnewline
\hline
\end{tabular}\end{center}
\end{table}

\renewcommand{\baselinestretch}{1}\normalsize
\begin{table}\begin{center}
\caption{Simulated size/power rates for the nominal 5\% level of Cramer -- von Misses tests of Models I-IV with static and dynamic specifications applied to simulated data, $T=100$.}
\label{t:rejcvm100}
\footnotesize
\renewcommand{\arraystretch}{0.5}
\begin{tabular}{cccccccccccc}\hline\hline
 &  &  & $\hat S_{2T} $  & $\hat R_{2T,50} $  & $\hat R_{2T,25} $  & $\hat R_{2T} $  & $\hat S_{1T} $  & $\hat R_{1T,50} $  & $\hat R_{1T,25} $  & $\hat R_{1T} $  & $\hat Z_{T} $ \tabularnewline
\hline\multicolumn{11}{c}{Size 1  $H_0:$ static probit} \tabularnewline
 & Model I &     & $ \tttt{5.5} $  & $ \tttt{6.0} $  & $ \tttt{5.5} $  & $ \tttt{4.5} $  & $ \tttt{6.6} $  & $ \tttt{6.3} $  & $ \tttt{5.7} $  & $ \tttt{5.4} $  & $ \tttt{7.8} $ \tabularnewline
 & Model II &     & $ \tttt{5.3} $  & $ \tttt{6.7} $  & $ \tttt{5.0} $  & $ \tttt{5.5} $  & $ \tttt{6.3} $  & $ \tttt{5.2} $  & $ \tttt{4.2} $  & $ \tttt{3.3} $  & $ \tttt{6.5} $ \tabularnewline
 & Model III &     & $ \tttt{7.7} $  & $ \tttt{7.0} $  & $ \tttt{6.5} $  & $ \tttt{5.4} $  & $ \tttt{6.0} $  & $ \tttt{3.7} $  & $ \tttt{3.3} $  & $ \tttt{4.5} $  & $ \tttt{6.4} $ \tabularnewline
 & Model IV &     & $ \tttt{5.2} $  & $ \tttt{6.7} $  & $ \tttt{5.6} $  & $ \tttt{3.9} $  & $ \tttt{5.1} $  & $ \tttt{4.6} $  & $ \tttt{4.9} $  & $ \tttt{2.8} $  & $ \tttt{6.4} $ \tabularnewline
\hline\multicolumn{11}{c}{Size 2  $H_0:$ static logit} \tabularnewline
 & Model I &     & $ \tttt{6.5} $  & $ \tttt{6.5} $  & $ \tttt{4.9} $  & $ \tttt{4.1} $  & $ \tttt{7.2} $  & $ \tttt{5.6} $  & $ \tttt{6.0} $  & $ \tttt{4.4} $  & $ \tttt{7.2} $ \tabularnewline
 & Model II &     & $ \tttt{5.6} $  & $ \tttt{6.7} $  & $ \tttt{7.6} $  & $ \tttt{4.0} $  & $ \tttt{4.6} $  & $ \tttt{5.3} $  & $ \tttt{4.6} $  & $ \tttt{4.8} $  & $ \tttt{5.6} $ \tabularnewline
 & Model III &     & $ \tttt{7.3} $  & $ \tttt{9.0} $  & $ \tttt{6.4} $  & $ \tttt{3.3} $  & $ \tttt{6.4} $  & $ \tttt{7.8} $  & $ \tttt{5.2} $  & $ \tttt{3.3} $  & $ \tttt{8.5} $ \tabularnewline
 & Model IV &     & $ \tttt{6.6} $  & $ \tttt{6.3} $  & $ \tttt{5.0} $  & $ \tttt{4.5} $  & $ \tttt{6.5} $  & $ \tttt{4.6} $  & $ \tttt{4.7} $  & $ \tttt{4.7} $  & $ \tttt{9.1} $ \tabularnewline
\hline\multicolumn{11}{c}{Power 1  $H_0:$ static probit vs  $H_1:$ static logit} \tabularnewline
 & Model I &     & $ \tttt{8.5} $  & $ \tttt{7.7} $  & $ \tttt{6.6} $  & $ \tttt{4.9} $  & $ \tttt{8.4} $  & $ \tttt{6.5} $  & $ \tttt{6.0} $  & $ \tttt{3.6} $  & $ \tttt{7.1} $ \tabularnewline
 & Model II &     & $ \tttt{5.1} $  & $ \tttt{5.0} $  & $ \tttt{4.4} $  & $ \tttt{4.0} $  & $ \tttt{6.4} $  & $ \tttt{6.9} $  & $ \tttt{5.3} $  & $ \tttt{4.0} $  & $ \tttt{8.7} $ \tabularnewline
 & Model III &     & $ \tttt{9.1} $  & $ \tttt{9.4} $  & $ \tttt{7.9} $  & $ \tttt{4.7} $  & $ \tttt{9.0} $  & $ \tttt{8.3} $  & $ \tttt{7.7} $  & $ \tttt{4.6} $  & $ \tttt{8.2} $ \tabularnewline
 & Model IV &     & $ \tttt{6.3} $  & $ \tttt{6.2} $  & $ \tttt{5.3} $  & $ \tttt{4.5} $  & $ \tttt{10.2} $  & $ \tttt{8.6} $  & $ \tttt{7.5} $  & $ \tttt{3.8} $  & $ \tttt{8.3} $ \tabularnewline
\hline\multicolumn{11}{c}{Power 2  $H_0:$ static probit vs  $H_1:$ dynamic probit} \tabularnewline
 & Model I &     & $ \tttt{89.2} $  & $ \tttt{85.2} $  & $ \tttt{82.7} $  & $ \tttt{25.7} $  & $ \tttt{13.2} $  & $ \tttt{12.0} $  & $ \tttt{11.7} $  & $ \tttt{4.6} $  & $ \tttt{18.4} $ \tabularnewline
 & Model II &     & $ \tttt{92.8} $  & $ \tttt{92.3} $  & $ \tttt{91.1} $  & $ \tttt{34.2} $  & $ \tttt{10.5} $  & $ \tttt{8.1} $  & $ \tttt{8.8} $  & $ \tttt{3.0} $  & $ \tttt{17.2} $ \tabularnewline
 & Model III &     & $ \tttt{90.7} $  & $ \tttt{88.4} $  & $ \tttt{86.1} $  & $ \tttt{22.5} $  & $ \tttt{9.2} $  & $ \tttt{9.8} $  & $ \tttt{8.5} $  & $ \tttt{5.0} $  & $ \tttt{9.4} $ \tabularnewline
 & Model IV &     & $ \tttt{88.1} $  & $ \tttt{84.1} $  & $ \tttt{83.0} $  & $ \tttt{27.7} $  & $ \tttt{10.3} $  & $ \tttt{7.8} $  & $ \tttt{7.4} $  & $ \tttt{4.4} $  & $ \tttt{12.5} $ \tabularnewline
\hline\multicolumn{11}{c}{Power 3  $H_0:$ static probit vs  $H_1:$ dynamic logit} \tabularnewline
 & Model I &     & $ \tttt{90.1} $  & $ \tttt{89.3} $  & $ \tttt{86.0} $  & $ \tttt{22.9} $  & $ \tttt{12.1} $  & $ \tttt{10.0} $  & $ \tttt{8.5} $  & $ \tttt{5.0} $  & $ \tttt{12.6} $ \tabularnewline
 & Model II &     & $ \tttt{94.2} $  & $ \tttt{93.0} $  & $ \tttt{90.6} $  & $ \tttt{29.8} $  & $ \tttt{9.6} $  & $ \tttt{9.1} $  & $ \tttt{7.1} $  & $ \tttt{3.9} $  & $ \tttt{14.6} $ \tabularnewline
 & Model III &     & $ \tttt{93.5} $  & $ \tttt{91.9} $  & $ \tttt{90.9} $  & $ \tttt{30.3} $  & $ \tttt{10.0} $  & $ \tttt{8.0} $  & $ \tttt{7.8} $  & $ \tttt{4.4} $  & $ \tttt{10.9} $ \tabularnewline
 & Model IV &     & $ \tttt{91.1} $  & $ \tttt{88.4} $  & $ \tttt{85.9} $  & $ \tttt{26.0} $  & $ \tttt{11.1} $  & $ \tttt{12.3} $  & $ \tttt{11.4} $  & $ \tttt{4.7} $  & $ \tttt{14.7} $ \tabularnewline
\hline
\end{tabular}
\end{center}\end{table}

\begin{table}\begin{center}
\caption{Simulated size/power  rates for the nominal 5\% level of  Kolmogorov -- Smirnov tests of Models I-IV with static and dynamic specifications applied to simulated data, $T=100$.}
\label{t:rejks100}
\footnotesize
\renewcommand{\arraystretch}{0.5}
\begin{tabular}{cccccccccccc}
\hline
\hline
 &  &  & $\hat S_{2T} $  & $\hat R_{2T,50} $  & $\hat R_{2T,25} $  & $\hat R_{2T} $  & $\hat S_{1T} $  & $\hat R_{1T,50} $  & $\hat R_{1T,25} $  & $\hat R_{1T} $  & $\hat Z_{T} $ \tabularnewline\hline\multicolumn{11}{c}{Size 1  $H_0:$ static probit} \tabularnewline
 & Model I &     & $ \tttt{5.1} $  & $ \tttt{6.4} $  & $ \tttt{5.2} $  & $ \tttt{3.9} $  & $ \tttt{7.8} $  & $ \tttt{6.3} $  & $ \tttt{6.8} $  & $ \tttt{4.9} $  & $ \tttt{7.9} $ \tabularnewline
 & Model II &     & $ \tttt{5.5} $  & $ \tttt{6.5} $  & $ \tttt{3.9} $  & $ \tttt{4.9} $  & $ \tttt{5.9} $  & $ \tttt{5.1} $  & $ \tttt{4.1} $  & $ \tttt{4.8} $  & $ \tttt{6.2} $ \tabularnewline
 & Model III &     & $ \tttt{7.7} $  & $ \tttt{7.8} $  & $ \tttt{6.8} $  & $ \tttt{5.1} $  & $ \tttt{6.1} $  & $ \tttt{7.0} $  & $ \tttt{6.0} $  & $ \tttt{4.9} $  & $ \tttt{5.6} $ \tabularnewline
 & Model IV &     & $ \tttt{6.5} $  & $ \tttt{5.4} $  & $ \tttt{5.3} $  & $ \tttt{3.4} $  & $ \tttt{5.3} $  & $ \tttt{5.3} $  & $ \tttt{4.8} $  & $ \tttt{3.6} $  & $ \tttt{7.2} $ \tabularnewline
\hline\multicolumn{11}{c}{Size 2  $H_0:$ static logit} \tabularnewline
 & Model I &     & $ \tttt{7.0} $  & $ \tttt{6.4} $  & $ \tttt{6.1} $  & $ \tttt{5.4} $  & $ \tttt{9.1} $  & $ \tttt{6.4} $  & $ \tttt{6.3} $  & $ \tttt{3.7} $  & $ \tttt{6.7} $ \tabularnewline
 & Model II &     & $ \tttt{4.7} $  & $ \tttt{4.9} $  & $ \tttt{4.6} $  & $ \tttt{3.5} $  & $ \tttt{5.6} $  & $ \tttt{3.8} $  & $ \tttt{4.0} $  & $ \tttt{4.8} $  & $ \tttt{5.8} $ \tabularnewline
 & Model III &     & $ \tttt{8.3} $  & $ \tttt{8.3} $  & $ \tttt{6.7} $  & $ \tttt{3.2} $  & $ \tttt{6.2} $  & $ \tttt{5.7} $  & $ \tttt{3.5} $  & $ \tttt{4.0} $  & $ \tttt{10.0} $ \tabularnewline
 & Model IV &     & $ \tttt{6.2} $  & $ \tttt{6.5} $  & $ \tttt{5.1} $  & $ \tttt{4.7} $  & $ \tttt{6.6} $  & $ \tttt{5.8} $  & $ \tttt{5.3} $  & $ \tttt{4.0} $  & $ \tttt{8.1} $ \tabularnewline
\hline\multicolumn{11}{c}{Power 1  $H_0:$ static probit vs  $H_1:$ static logit} \tabularnewline
 & Model I &     & $ \tttt{7.0} $  & $ \tttt{6.2} $  & $ \tttt{5.4} $  & $ \tttt{3.7} $  & $ \tttt{5.2} $  & $ \tttt{3.3} $  & $ \tttt{3.9} $  & $ \tttt{3.2} $  & $ \tttt{7.7} $ \tabularnewline
 & Model II &     & $ \tttt{4.3} $  & $ \tttt{3.8} $  & $ \tttt{4.5} $  & $ \tttt{3.7} $  & $ \tttt{4.1} $  & $ \tttt{3.9} $  & $ \tttt{3.6} $  & $ \tttt{3.9} $  & $ \tttt{8.9} $ \tabularnewline
 & Model III &     & $ \tttt{10.2} $  & $ \tttt{7.3} $  & $ \tttt{7.1} $  & $ \tttt{3.9} $  & $ \tttt{7.1} $  & $ \tttt{5.7} $  & $ \tttt{5.7} $  & $ \tttt{4.5} $  & $ \tttt{9.2} $ \tabularnewline
 & Model IV &     & $ \tttt{5.6} $  & $ \tttt{6.6} $  & $ \tttt{4.3} $  & $ \tttt{3.2} $  & $ \tttt{6.4} $  & $ \tttt{5.1} $  & $ \tttt{6.2} $  & $ \tttt{3.4} $  & $ \tttt{6.8} $ \tabularnewline
\hline\multicolumn{11}{c}{Power 2  $H_0:$ static probit vs  $H_1:$ dynamic probit} \tabularnewline
 & Model I &     & $ \tttt{82.8} $  & $ \tttt{79.0} $  & $ \tttt{74.5} $  & $ \tttt{13.6} $  & $ \tttt{10.3} $  & $ \tttt{9.1} $  & $ \tttt{7.1} $  & $ \tttt{3.5} $  & $ \tttt{16.9} $ \tabularnewline
 & Model II &     & $ \tttt{87.9} $  & $ \tttt{85.5} $  & $ \tttt{83.3} $  & $ \tttt{17.7} $  & $ \tttt{12.1} $  & $ \tttt{11.2} $  & $ \tttt{9.3} $  & $ \tttt{3.3} $  & $ \tttt{14.0} $ \tabularnewline
 & Model III &     & $ \tttt{85.7} $  & $ \tttt{83.2} $  & $ \tttt{79.4} $  & $ \tttt{13.8} $  & $ \tttt{7.1} $  & $ \tttt{6.4} $  & $ \tttt{7.2} $  & $ \tttt{3.9} $  & $ \tttt{9.2} $ \tabularnewline
 & Model IV &     & $ \tttt{81.7} $  & $ \tttt{78.5} $  & $ \tttt{74.6} $  & $ \tttt{13.8} $  & $ \tttt{7.7} $  & $ \tttt{7.8} $  & $ \tttt{6.7} $  & $ \tttt{4.9} $  & $ \tttt{11.3} $ \tabularnewline
\hline\multicolumn{11}{c}{Power 3  $H_0:$ static probit vs  $H_1:$ dynamic logit} \tabularnewline
 & Model I &     & $ \tttt{86.2} $  & $ \tttt{82.7} $  & $ \tttt{79.0} $  & $ \tttt{14.2} $  & $ \tttt{7.7} $  & $ \tttt{4.9} $  & $ \tttt{3.8} $  & $ \tttt{4.2} $  & $ \tttt{11.8} $ \tabularnewline
 & Model II &     & $ \tttt{90.0} $  & $ \tttt{86.2} $  & $ \tttt{82.2} $  & $ \tttt{15.9} $  & $ \tttt{9.3} $  & $ \tttt{7.9} $  & $ \tttt{8.1} $  & $ \tttt{4.1} $  & $ \tttt{14.2} $ \tabularnewline
 & Model III &     & $ \tttt{89.0} $  & $ \tttt{86.4} $  & $ \tttt{83.7} $  & $ \tttt{15.9} $  & $ \tttt{5.6} $  & $ \tttt{5.1} $  & $ \tttt{4.4} $  & $ \tttt{4.6} $  & $ \tttt{10.5} $ \tabularnewline
 & Model IV &     & $ \tttt{87.5} $  & $ \tttt{83.8} $  & $ \tttt{79.3} $  & $ \tttt{16.1} $  & $ \tttt{9.4} $  & $ \tttt{7.5} $  & $ \tttt{7.7} $  & $ \tttt{5.9} $  & $ \tttt{12.9} $ \tabularnewline
\hline
\end{tabular}
\end{center}\end{table}


\begin{table}\begin{center}
\caption{Simulated size/power rates for the nominal 5\% level of  Cramer -- von Misses tests of Models I-IV with static and dynamic specifications applied to simulated data, $T=200$.}
\label{t:rejcvm200}
\footnotesize
\renewcommand{\arraystretch}{0.5}
\begin{tabular}{cccccccccccc}
\hline
\hline
 &  &  & $\hat S_{2T} $  & $\hat R_{2T,50} $  & $\hat R_{2T,25} $  & $\hat R_{2T} $  & $\hat S_{1T} $  & $\hat R_{1T,50} $  & $\hat R_{1T,25} $  & $\hat R_{1T} $  & $\hat Z_{T} $ \tabularnewline\hline\multicolumn{11}{c}{Size 1  $H_0:$ static probit} \tabularnewline
 & Model I &     & $ \tttt{4.0} $  & $ \tttt{5.4} $  & $ \tttt{5.7} $  & $ \tttt{6.2} $  & $ \tttt{4.2} $  & $ \tttt{4.6} $  & $ \tttt{4.9} $  & $ \tttt{5.8} $  & $ \tttt{5.2} $ \tabularnewline
 & Model II &     & $ \tttt{4.5} $  & $ \tttt{4.4} $  & $ \tttt{3.5} $  & $ \tttt{2.4} $  & $ \tttt{6.3} $  & $ \tttt{4.7} $  & $ \tttt{5.9} $  & $ \tttt{4.4} $  & $ \tttt{7.0} $ \tabularnewline
 & Model III &     & $ \tttt{4.6} $  & $ \tttt{4.4} $  & $ \tttt{3.4} $  & $ \tttt{4.2} $  & $ \tttt{5.4} $  & $ \tttt{5.5} $  & $ \tttt{5.2} $  & $ \tttt{3.3} $  & $ \tttt{5.4} $ \tabularnewline
 & Model IV &     & $ \tttt{5.3} $  & $ \tttt{6.1} $  & $ \tttt{6.3} $  & $ \tttt{4.4} $  & $ \tttt{4.8} $  & $ \tttt{4.6} $  & $ \tttt{7.0} $  & $ \tttt{4.9} $  & $ \tttt{6.9} $ \tabularnewline
\hline\multicolumn{11}{c}{Size 2  $H_0:$ static logit} \tabularnewline
 & Model I &     & $ \tttt{7.2} $  & $ \tttt{8.2} $  & $ \tttt{6.7} $  & $ \tttt{5.8} $  & $ \tttt{5.8} $  & $ \tttt{6.7} $  & $ \tttt{6.4} $  & $ \tttt{3.8} $  & $ \tttt{5.2} $ \tabularnewline
 & Model II &     & $ \tttt{5.4} $  & $ \tttt{6.4} $  & $ \tttt{6.1} $  & $ \tttt{4.7} $  & $ \tttt{4.8} $  & $ \tttt{5.6} $  & $ \tttt{5.3} $  & $ \tttt{5.2} $  & $ \tttt{6.1} $ \tabularnewline
 & Model III &     & $ \tttt{5.3} $  & $ \tttt{5.2} $  & $ \tttt{3.9} $  & $ \tttt{4.0} $  & $ \tttt{5.6} $  & $ \tttt{5.8} $  & $ \tttt{6.7} $  & $ \tttt{4.4} $  & $ \tttt{6.9} $ \tabularnewline
 & Model IV &     & $ \tttt{5.4} $  & $ \tttt{6.8} $  & $ \tttt{5.0} $  & $ \tttt{4.0} $  & $ \tttt{5.6} $  & $ \tttt{5.2} $  & $ \tttt{5.0} $  & $ \tttt{4.1} $  & $ \tttt{8.3} $ \tabularnewline
\hline\multicolumn{11}{c}{Power 1  $H_0:$ static probit vs  $H_1:$ static logit} \tabularnewline
 & Model I &     & $ \tttt{7.2} $  & $ \tttt{8.2} $  & $ \tttt{6.6} $  & $ \tttt{6.9} $  & $ \tttt{10.9} $  & $ \tttt{10.3} $  & $ \tttt{10.9} $  & $ \tttt{6.9} $  & $ \tttt{9.2} $ \tabularnewline
 & Model II &     & $ \tttt{4.5} $  & $ \tttt{4.9} $  & $ \tttt{5.6} $  & $ \tttt{6.3} $  & $ \tttt{7.5} $  & $ \tttt{6.4} $  & $ \tttt{7.3} $  & $ \tttt{6.7} $  & $ \tttt{6.5} $ \tabularnewline
 & Model III &     & $ \tttt{6.0} $  & $ \tttt{5.2} $  & $ \tttt{6.1} $  & $ \tttt{5.9} $  & $ \tttt{6.9} $  & $ \tttt{6.8} $  & $ \tttt{7.9} $  & $ \tttt{6.6} $  & $ \tttt{7.0} $ \tabularnewline
 & Model IV &     & $ \tttt{6.5} $  & $ \tttt{6.6} $  & $ \tttt{6.6} $  & $ \tttt{4.3} $  & $ \tttt{7.1} $  & $ \tttt{6.1} $  & $ \tttt{7.5} $  & $ \tttt{5.8} $  & $ \tttt{5.4} $ \tabularnewline
\hline\multicolumn{11}{c}{Power 2  $H_0:$ static probit vs  $H_1:$ dynamic probit} \tabularnewline
 & Model I &     & $ \tttt{98.5} $  & $ \tttt{97.3} $  & $ \tttt{95.2} $  & $ \tttt{33.2} $  & $ \tttt{13.6} $  & $ \tttt{11.6} $  & $ \tttt{9.8} $  & $ \tttt{7.5} $  & $ \tttt{16.2} $ \tabularnewline
 & Model II &     & $ \tttt{99.5} $  & $ \tttt{99.3} $  & $ \tttt{98.5} $  & $ \tttt{41.5} $  & $ \tttt{16.0} $  & $ \tttt{14.8} $  & $ \tttt{12.6} $  & $ \tttt{7.1} $  & $ \tttt{18.2} $ \tabularnewline
 & Model III &     & $ \tttt{98.5} $  & $ \tttt{97.0} $  & $ \tttt{95.9} $  & $ \tttt{30.7} $  & $ \tttt{13.0} $  & $ \tttt{11.8} $  & $ \tttt{9.8} $  & $ \tttt{7.9} $  & $ \tttt{13.9} $ \tabularnewline
 & Model IV &     & $ \tttt{95.8} $  & $ \tttt{93.6} $  & $ \tttt{91.6} $  & $ \tttt{22.9} $  & $ \tttt{10.2} $  & $ \tttt{9.5} $  & $ \tttt{7.6} $  & $ \tttt{5.3} $  & $ \tttt{13.7} $ \tabularnewline
\hline\multicolumn{11}{c}{Power 3  $H_0:$ static probit vs  $H_1:$ dynamic logit} \tabularnewline
 & Model I &     & $ \tttt{98.6} $  & $ \tttt{97.5} $  & $ \tttt{95.6} $  & $ \tttt{34.5} $  & $ \tttt{15.2} $  & $ \tttt{14.0} $  & $ \tttt{14.0} $  & $ \tttt{5.4} $  & $ \tttt{16.7} $ \tabularnewline
 & Model II &     & $ \tttt{99.5} $  & $ \tttt{98.9} $  & $ \tttt{98.6} $  & $ \tttt{39.1} $  & $ \tttt{16.9} $  & $ \tttt{16.1} $  & $ \tttt{13.7} $  & $ \tttt{7.4} $  & $ \tttt{20.8} $ \tabularnewline
 & Model III &     & $ \tttt{98.8} $  & $ \tttt{98.1} $  & $ \tttt{96.4} $  & $ \tttt{31.2} $  & $ \tttt{14.8} $  & $ \tttt{13.7} $  & $ \tttt{11.6} $  & $ \tttt{6.6} $  & $ \tttt{17.9} $ \tabularnewline
 & Model IV &     & $ \tttt{95.8} $  & $ \tttt{94.2} $  & $ \tttt{91.8} $  & $ \tttt{23.9} $  & $ \tttt{11.4} $  & $ \tttt{10.1} $  & $ \tttt{8.6} $  & $ \tttt{5.3} $  & $ \tttt{11.4} $ \tabularnewline
\hline
\end{tabular}\end{center}
\end{table}

\begin{table}\begin{center}
\caption{Simulated size/power rates for the nominal 5\% level of  Kolmogorov -- Smirnov tests of Models I-IV with static and dynamic specifications applied to simulated data, $T=200$.}
\label{t:rejks200}
\footnotesize
\renewcommand{\arraystretch}{0.5}
\begin{tabular}{cccccccccccc}
\hline
\hline
 &  &  & $\hat S_{2T} $  & $\hat R_{2T,50} $  & $\hat R_{2T,25} $  & $\hat R_{2T} $  & $\hat S_{1T} $  & $\hat R_{1T,50} $  & $\hat R_{1T,25} $  & $\hat R_{1T} $  & $\hat Z_{T} $ \tabularnewline\hline\multicolumn{11}{c}{Size 1  $H_0:$ static probit} \tabularnewline
 & Model I &     & $ \tttt{5.1} $  & $ \tttt{5.0} $  & $ \tttt{6.0} $  & $ \tttt{4.8} $  & $ \tttt{4.5} $  & $ \tttt{5.9} $  & $ \tttt{3.9} $  & $ \tttt{5.1} $  & $ \tttt{5.5} $ \tabularnewline
 & Model II &     & $ \tttt{3.7} $  & $ \tttt{3.9} $  & $ \tttt{3.9} $  & $ \tttt{2.9} $  & $ \tttt{6.3} $  & $ \tttt{5.6} $  & $ \tttt{6.4} $  & $ \tttt{4.6} $  & $ \tttt{5.7} $ \tabularnewline
 & Model III &     & $ \tttt{4.5} $  & $ \tttt{5.2} $  & $ \tttt{4.3} $  & $ \tttt{3.9} $  & $ \tttt{4.2} $  & $ \tttt{4.5} $  & $ \tttt{4.9} $  & $ \tttt{4.3} $  & $ \tttt{5.3} $ \tabularnewline
 & Model IV &     & $ \tttt{5.0} $  & $ \tttt{6.4} $  & $ \tttt{7.0} $  & $ \tttt{4.6} $  & $ \tttt{4.6} $  & $ \tttt{4.8} $  & $ \tttt{6.4} $  & $ \tttt{6.9} $  & $ \tttt{6.8} $ \tabularnewline
\hline\multicolumn{11}{c}{Size 2  $H_0:$ static logit} \tabularnewline
 & Model I &     & $ \tttt{5.7} $  & $ \tttt{5.7} $  & $ \tttt{6.3} $  & $ \tttt{4.9} $  & $ \tttt{6.3} $  & $ \tttt{5.9} $  & $ \tttt{6.3} $  & $ \tttt{3.5} $  & $ \tttt{4.8} $ \tabularnewline
 & Model II &     & $ \tttt{5.5} $  & $ \tttt{5.1} $  & $ \tttt{5.9} $  & $ \tttt{3.4} $  & $ \tttt{5.4} $  & $ \tttt{4.6} $  & $ \tttt{5.9} $  & $ \tttt{4.9} $  & $ \tttt{5.3} $ \tabularnewline
 & Model III &     & $ \tttt{3.6} $  & $ \tttt{5.4} $  & $ \tttt{4.6} $  & $ \tttt{3.6} $  & $ \tttt{6.4} $  & $ \tttt{4.3} $  & $ \tttt{5.2} $  & $ \tttt{5.3} $  & $ \tttt{7.8} $ \tabularnewline
 & Model IV &     & $ \tttt{6.4} $  & $ \tttt{7.3} $  & $ \tttt{5.6} $  & $ \tttt{4.7} $  & $ \tttt{6.6} $  & $ \tttt{6.4} $  & $ \tttt{4.7} $  & $ \tttt{4.5} $  & $ \tttt{8.5} $ \tabularnewline
\hline\multicolumn{11}{c}{Power 1  $H_0:$ static probit vs  $H_1:$ static logit} \tabularnewline
 & Model I &     & $ \tttt{6.3} $  & $ \tttt{6.5} $  & $ \tttt{4.2} $  & $ \tttt{6.6} $  & $ \tttt{7.2} $  & $ \tttt{5.9} $  & $ \tttt{5.0} $  & $ \tttt{6.7} $  & $ \tttt{8.7} $ \tabularnewline
 & Model II &     & $ \tttt{4.6} $  & $ \tttt{5.0} $  & $ \tttt{6.4} $  & $ \tttt{6.3} $  & $ \tttt{5.2} $  & $ \tttt{4.9} $  & $ \tttt{5.7} $  & $ \tttt{6.5} $  & $ \tttt{6.1} $ \tabularnewline
 & Model III &     & $ \tttt{5.0} $  & $ \tttt{6.2} $  & $ \tttt{5.7} $  & $ \tttt{5.2} $  & $ \tttt{3.7} $  & $ \tttt{4.1} $  & $ \tttt{5.0} $  & $ \tttt{6.3} $  & $ \tttt{7.1} $ \tabularnewline
 & Model IV &     & $ \tttt{5.7} $  & $ \tttt{7.0} $  & $ \tttt{5.6} $  & $ \tttt{4.5} $  & $ \tttt{5.4} $  & $ \tttt{3.4} $  & $ \tttt{4.6} $  & $ \tttt{6.1} $  & $ \tttt{5.1} $ \tabularnewline
\hline\multicolumn{11}{c}{Power 2  $H_0:$ static probit vs  $H_1:$ dynamic probit} \tabularnewline
 & Model I &     & $ \tttt{94.3} $  & $ \tttt{92.0} $  & $ \tttt{86.5} $  & $ \tttt{22.8} $  & $ \tttt{11.4} $  & $ \tttt{10.6} $  & $ \tttt{9.7} $  & $ \tttt{5.9} $  & $ \tttt{14.0} $ \tabularnewline
 & Model II &     & $ \tttt{98.1} $  & $ \tttt{96.5} $  & $ \tttt{94.4} $  & $ \tttt{26.3} $  & $ \tttt{15.5} $  & $ \tttt{13.1} $  & $ \tttt{13.5} $  & $ \tttt{7.3} $  & $ \tttt{13.1} $ \tabularnewline
 & Model III &     & $ \tttt{94.3} $  & $ \tttt{91.0} $  & $ \tttt{87.9} $  & $ \tttt{21.0} $  & $ \tttt{14.7} $  & $ \tttt{13.1} $  & $ \tttt{12.7} $  & $ \tttt{7.3} $  & $ \tttt{13.8} $ \tabularnewline
 & Model IV &     & $ \tttt{90.5} $  & $ \tttt{85.0} $  & $ \tttt{82.0} $  & $ \tttt{17.9} $  & $ \tttt{11.0} $  & $ \tttt{9.5} $  & $ \tttt{9.4} $  & $ \tttt{6.3} $  & $ \tttt{11.4} $ \tabularnewline
\hline\multicolumn{11}{c}{Power 3  $H_0:$ static probit vs  $H_1:$ dynamic logit} \tabularnewline
 & Model I &     & $ \tttt{97.1} $  & $ \tttt{93.8} $  & $ \tttt{91.8} $  & $ \tttt{24.7} $  & $ \tttt{12.4} $  & $ \tttt{12.8} $  & $ \tttt{11.1} $  & $ \tttt{5.5} $  & $ \tttt{13.4} $ \tabularnewline
 & Model II &     & $ \tttt{98.9} $  & $ \tttt{97.6} $  & $ \tttt{96.5} $  & $ \tttt{29.5} $  & $ \tttt{16.9} $  & $ \tttt{17.1} $  & $ \tttt{14.6} $  & $ \tttt{7.4} $  & $ \tttt{16.9} $ \tabularnewline
 & Model III &     & $ \tttt{96.1} $  & $ \tttt{93.8} $  & $ \tttt{91.8} $  & $ \tttt{26.0} $  & $ \tttt{14.6} $  & $ \tttt{14.4} $  & $ \tttt{11.9} $  & $ \tttt{8.0} $  & $ \tttt{15.4} $ \tabularnewline
 & Model IV &     & $ \tttt{93.0} $  & $ \tttt{89.2} $  & $ \tttt{86.9} $  & $ \tttt{14.1} $  & $ \tttt{13.0} $  & $ \tttt{12.5} $  & $ \tttt{8.5} $  & $ \tttt{5.8} $  & $ \tttt{10.4} $ \tabularnewline
\hline
\end{tabular}\end{center}
\end{table}
\renewcommand{\baselinestretch}{1.5}\normalsize

Parameter estimates for real data are reported in Tables~\ref{t:estprobit}
and~\ref{t:estlogit}. The main question is whether the static Probit or
Logit models are appropriate for changes in the interest rates, and we check
this with our tests. The $p$-values in Tables~\ref{t:pvalcvm} and~\ref%
{t:pvalks} say that all these models are rejected even at the 1\%
significance level by biparameter nonrandomized transform based tests. Note
that single parameter static tests (e.g. $\hat R_{1T}$, $\hat S_{1T}$)
cannot reject any proposed model with the sole exception of $\hat S_{1T}$
which rejects at 5\% Model II with Cramer -- von Misses test statistics.

To study the reliability of these results we conduct a Monte Carlo
experiment using the estimated models with the real data as data generating
processes and obtain the simulations for the discrete response conditional
on the covariates time series. In Tables~\ref{t:rejcvm100} and~\ref%
{t:rejks100} we provide the empirical size and power results of our tests
across simulations for sample size $T=100$ and static Probit and Logit and
output gap choices (Models I to IV). To speed up the
simulation procedure, we use the warp bootstrap algorithm of Giacomini,
Politis and White (2013). We see that all bootstrap tests provide
reasonable size accuracy, tests based on single parameter empirical
processes underrejecting slightly, while ones based on bivariate processes
tend to overreject moderately. Kolmogorov-Smirnov and Cramer-von Mises tests
perform similarly in all cases, and the choice of the output gap series does
not make large differences either, nor does the introduction of lagged
endogenous (discrete) variables in the information set.

The power of the tests for the static Probit model is analyzed against three
different alternatives: static Logit, dynamic Probit and dynamic Logit. We
see that the tests without randomization, $\hat{S}_{1T}$ and $\hat{S}_{2T}$
always perform better than random continuous processes ${\hat{R}_{1T,M}}$
and ${\hat{R}}_{2T,M}$, which in turn dominate ${\hat{R}_{1T}}$ and ${\hat{R}%
}_{2T}$, thus confirming our theoretical findings. When we compare Probit
and Logit specifications while letting the dynamic aspect of the model be
well specified, static in both cases, we observe that with this sample size
and these specifications, it is almost impossible to distinguish Probit from
Logit models. The power against a dynamic Probit and Logit alternatives is
very high. Since the nature of misspecification is dynamic, once again
bivariate processes should have more power compared to single parameter
counterparts, as it is confirmed in our simulation results. It can also be
observed that for these alternatives, the Cramer-von Mises criterium
provides more power than Kolmogorov-Smirnov tests. As for alternative tests
based on $\hat{Z}_{T}$, they have power comparable to $\hat{S}_{1T}$,
sometimes slightly better, and are always outperformed by any bivariate
test. This is not surprising, since $\hat{Z}_{T}$ has more structure, i.e.
it assumes a single-index model for covariates, but averages across points,
thus suffering the same problems as other single parameter tests considered
here.

In Tables~\ref{t:rejcvm200} and~\ref{t:rejks200} we provide the empirical
size and power results of our tests for the larger sample size $T=200$. Here
the size properties are similar, while power rejections rates are noticeably
higher for the dynamic alternatives.

\section{CONCLUSIONS}

In this paper we have proposed new specification tests for the conditional
distribution of discrete data with possibly infinite support. The new tests are functionals of
empirical processes based on a nonrandomized transform that solves the
implementation problem of the usual PIT for discrete distributions and
achieves consistency against a wide class of alternatives. We show the
validity of a bootstrap algorithm for approximating the null distribution of
the test statistics, which are model and parameter dependent. In our
simulation study, we show that our method compares favorably in many
relevant situations with other methods available in the literature and have
illustrated the new method in a small application.

\section{APPENDIX}

\subsection{Properties of the nonrandomized transform}

In this section we derive the basic properties of the nonrandomized
transform, which are required prior to proving the weak convergence results
for our empirical process. Without loss of generality and in order to make
the exposition more transparent, we omit subscripts $t,\theta _{0}$ and
conditioning set $\Omega _{t}$, and use shortcuts $I_{F}\left( Y,u\right)
=I_{t,\theta _{0}}\left( Y_{t},u\right) $ and $I_{F,M}\left( Y,u\right)
=I_{t,\theta _{0},M}\left( Y_{t},u\right) $.

For $F\in \mathop{\mathcal M}\nolimits$, $F\left( F^{-1}(u)\right) \geq
u>F\left( F^{-1}(u)-1\right) $ and equality holds iff $u=F(k)$ for some
integer $k$. For a random variable $Y\sim G\in \mathop{\mathcal M}\nolimits$
we find $\Pr_{G}\left( F\left( Y\right) <u\right) =G\left( F^{-1}\left(
u\right) -1\right) $ and $g\left( F^{-1}\left( u\right) \right)
:=\Pr_{G}\left( Y=F^{-1}\left( u\right) \right) =G\left( F^{-1}\left(
u\right) \right) -G\left( F^{-1}\left( u\right) -1\right) $. For $G=F$, we
have that $\Pr_{F}\left( F\left( Y\right) <u\right) =F\left( F^{-1}\left(
u\right) -1\right) <u$, i.e. $F\left( Y\right) $ is not uniform and the
expectation of the indicator function $I\left( F\left( Y\right) <u\right) $
is never $u$ as it is for continuous $F$.

The nonrandomized transform can be written as 
\begin{equation*}
I_{F}\left( Y,u\right) =\left( 1-\delta _{F}\left( u\right) \right) %
\mathbbm{1}\! \left \{ Y=F^{-1}\left( u\right) \right \} +\mathbbm{1}\!
\left \{ Y<F^{-1}\left( u\right) \right \} ,
\end{equation*}%
where 
\begin{equation*}
\delta _{F}\left( u\right) :=\frac{F\left( F^{-1}\left( u\right) \right) -u}{%
f\left( F^{-1}\left( u\right) \right) }.
\end{equation*}%
Note that $\delta _{F}\left( u\right) \in \lbrack 0,1)$. We see that $%
I_{F}\left( Y,u\right) $ is a piecewise linear (continuous) function
increasing in $u$.
\begin{table}\begin{center}
\caption{Values of functionals of the new nonrandomized transform $I\left(\cdot,\cdot\right)$ for all possible values of $Y$ relative to inverted cdfs at points $u$ and $v$. For instance, $I_F\left(Y,u\right)-I_F\left(Y,v\right)=0$ if $Y<F^{-1}\left(u\right)$ and $Y<F^{-1}\left(v\right)$, while $I_F\left(Y,u\right)-I_F\left(Y,v\right)=-\delta_F\left(u\right)$ if $Y=F^{-1}\left(u\right)<F^{-1}\left(v\right)$.}
\label{t:support}
\begin{tabular}{rccc}
\hline\hline
 &  $Y<F^{-1}\left(u\right)$ &    $Y=F^{-1}\left(u\right)$ &  $Y>F^{-1}\left(u\right)$ \\
\hline\multicolumn{4}{c}{The value of $I_F\left(Y,u\right)$} \tabularnewline
 &  $1$ &    $1-\delta_F\left(u\right)$ &  $0$ \\
\hline\multicolumn{4}{c}{The value of $\ind{I_F\left(Y,u\right)\le v}$} \tabularnewline
 $v=0$ & $0$  &    $0$ &  $1$ \\
 $v\in(0,1)$ &  $0$ &    $\ind{1-\delta_F\left(u\right)\le v}$ &  $1$ \\
 $v=1$ &  $1$ &    $1$ &  $1$ \\
\hline\multicolumn{4}{c}{The value of $I_F\left(Y,u\right)-I_F\left(Y,v\right)$}
\tabularnewline
 $Y<F^{-1}\left(v\right)$ & $0$  &    $-\delta_F\left(u\right)$ &  $-1$ \\
$Y=F^{-1}\left(v\right)$ &  $\delta_F\left(v\right)$ &    $\delta_F\left(v\right)-\delta_F\left(u\right)$ &  $-1+\delta_F\left(v\right)$ \\
$Y>F^{-1}\left(v\right)$ &  $1$ &    $1-\delta_F\left(u\right)$ &  $0$ \\
\hline\multicolumn{4}{c}{The value of $I_F\left(Y,u\right)I_F\left(Y,v\right)$}
\tabularnewline
 $Y<F^{-1}\left(v\right)$ & $1$  &    $1-\delta_F\left(u\right)$ &  $0$ \\
 $Y=F^{-1}\left(v\right)$ &  $1-\delta_F\left(v\right)$ &    $\left(1-\delta_F\left(u\right)\right)\left(1-\delta_F\left(v\right)\right)$ &  $0$ \\
 $Y>F^{-1}\left(v\right)$ &  $0$ &    $0$ &  $0$ \\
\hline\multicolumn{4}{c}{The value of $I_F\left(Y,u\right)-I_H\left(Y,u\right)$ }
\tabularnewline
$Y<H^{-1}\left(u\right)$ & $0$  &    $-\delta_F\left(u\right)$ &  $-1$ \\
 $Y=H^{-1}\left(u\right)$ &  $\delta_H\left(u\right)$ &    $\delta_H\left(u\right)-\delta_F\left(u\right)$ &  $-1+\delta_H\left(u\right)$ \\
$Y>H^{-1}\left(u\right)$ &  $1$ &    $1-\delta_F\left(u\right)$ &  $0$ \\
\hline\multicolumn{4}{c}{The value of $I_F\left(Y,u\right)I_H\left(Y,u\right)$ }
\tabularnewline
$Y<H^{-1}\left(u\right)$ & $1$  &    $1-\delta_F\left(u\right)$ &  $0$ \\
$Y=H^{-1}\left(u\right)$ &  $1-\delta_H\left(u\right)$ &    $\left(1-\delta_F\left(u\right)\right)\left(1-\delta_H\left(u\right)\right)$ &  $0$ \\
$Y>H^{-1}\left(u\right)$ &  $0$ &    $0$ &  $0$ \\ \hline
\end{tabular}\end{center}
\end{table}
Let 
\begin{equation*}
\begin{split}
\delta _{F}\left( u,v\right) :={}& \left( \delta _{F}\left( u\vee v\right)
-\delta _{F}\left( u\right) \delta _{F}\left( v\right) \right) f\left(
F^{-1}\left( u\wedge v\right) \right) \\
& \times \mathbbm{1}\! \left \{ F^{-1}\left( u\right) =F^{-1}\left( v\right)
\right \} \in \left[ 0,u\wedge v\wedge f\left( F^{-1}\left( u\wedge v\right)
\right) \right] ,
\end{split}%
\end{equation*}%
\begin{equation*}
\begin{split}
d\left( G,F,u,v\right) :={}& d\left( G,F,u\wedge v\right) \\
& -\left( \delta _{F}\left( u\vee v\right) -\delta _{F}\left( u\right)
\delta _{F}\left( v\right) \right) \mathbbm{1}\! \left \{ F^{-1}\left(
u\right) =F^{-1}\left( v\right) \right \} \\
& \quad \times \left( g\left( F^{-1}\left( u\right) \right) -f\left(
F^{-1}\left( u\right) \right) \right) ,
\end{split}%
\end{equation*}%

\bigskip

In Table~\ref{t:support} and Lemma~A we list the properties of this
transform.

\begin{lemmaa}
\label{lem:app} For $0\le v\le u\le 1$ and $F, G, H \in \mathop{\mathcal M}%
\nolimits,$
\begin{enumerate}[(i)]
\item \label{a:ei} $\mathop{\rm E}\nolimits_{G}\left[ I_{F}\left( Y,u\right) %
\right] =u+d\left( G,F,u\right) $, where $\mathop{\rm E}\nolimits_{G}\left[
\cdot \right] =\int (\cdot )dG$ and $d\left( G,F,u\right) \in \lbrack -u,1-u]
$. When $G=F$, the expectation is $u$.

\item \label{a:ii} $I_{F}\left( Y,u\right) I_{F}\left( Y,v\right)
=I_{F}\left( Y,u\wedge v\right) -$\newline
$\left( \delta _{F}\left( u\vee v\right) -\delta _{F}\left( u\right) \delta
_{F}\left( v\right) \right) \times \mathbbm{1}\!\left\{ Y=F^{-1}\left(
u\right) =F^{-1}\left( v\right) \right\} .$

\item \label{a:eii} $\mathop{\rm E}\nolimits_{G}\left[ I_{F}\left(
Y,u\right) I_{F}\left( Y,v\right) \right] =u\wedge v-\delta _{F}\left(
u,v\right) +d\left( G,F,u,v\right) $.

\item \label{a:ai-i} $\left\vert I_{F}\left( Y,u\right) -I_{H}\left(
Y,u\right) \right\vert \leq 1\wedge \frac{\left\vert F(Y)-H(Y)\right\vert
\vee \left\vert F(Y-1)-H(Y-1)\right\vert }{f(Y)\vee h(Y)}$\newline
Moreover, $\mathop{\rm E}\nolimits_{F}\left[ \left\vert I_{F}\left(
Y,u\right) -I_{H}\left( Y,u\right) \right\vert ^{2}\right] \leq 9\sup_{k}{%
\left\vert {F}\left( k\right) -{H}\left( k\right) \right\vert }$.

\item \label{a:i-i} $\left\vert I_{F}\left( Y,u\right) -u-I_{F}\left(
Y,v\right) +v\right\vert \leq |u-v|\vee (1-f(Y))\text{ and }\left\vert
I_{F}\left( Y,u\right) -u-I_{F}\left( Y,v\right) +v\right\vert =|u-v|\text{
if }u,v\leq F(Y-1)\text{ or }u,v\geq F(Y).$ \newline
Moreover, $\mathop{\rm E}\nolimits_{F}\left[ \sup_{u,v\in \Psi (\varepsilon
)}\left\vert I_{F}\left( Y,u\right) -u-I_{F}\left( Y,v\right) +v\right\vert
^{2}\right] \leq 4\varepsilon ^{2}$, for any interval $\Psi (\varepsilon
)\subset \lbrack 0,1]$ of length $\varepsilon ^{2}$.

\item \label{a:eid} $\mathop{\rm E}\nolimits_{F_{z}}\left[ \mathbbm{1}%
\!\left\{ F^{\dag }\left( Y^{\dag }\right) <u\right\} \right] =I_{F}\left(
Y,u\right) $.

\item \label{a:eiim} $\mathop{\rm E}\nolimits_{F_{z}}\left[ I_{F,M}\left(
Y,u\right) I_{F,M}\left( Y,v\right) \right] =\frac{1}{M}I_{F}\left(
Y,u\wedge v\right) +\left( 1-\frac{1}{M}\right) I_{F}\left( Y,u\right)
I_{F}\left( Y,v\right) $.
\end{enumerate}
\end{lemmaa}

\subsection{Functional weak convergence of discrete martingales}

In this section we present Lindeberg-Feller-type sufficient conditions for
functional weak convergence of discrete martingales. In general, to
establish the weak convergence one needs to check tightness and
finite-dimensional convergence. In case of martingales, both parts can be
verified without imposing restrictive conditions. Here we state a result of
Nishiyama (2000), which extends Theorem 2.11.9 of van der Vaart and Wellner
(1996) to martingales, see also Theorem A.1 in Delgado and Escanciano
(2007). 
Further details on notation and definitions can be found in books Van der
Vaart and Wellner (1996) for empirical processes and row-independent
triangular arrays and in Jacod and Shiryaev (2003) for finite-dimensional
semimartingales. For every $T$, let $\left(\Omega^T,\mathop{\mathcal F}%
\nolimits^T,\{ \mathop{\mathcal F}\nolimits^T_{t}\},P^T\right)$ be a
discrete stochastic basis, where $\left(\Omega^T,\mathop{\mathcal F}%
\nolimits^T,P^T\right)$ is a probability space equipped with a filtration $%
\left \{ \mathop{\mathcal F}\nolimits^T_{t}\right \}$. For nonempty set $\Psi
$, let $\{ \xi^T_{t}\}_{t=1,2,\ldots}$ be a $\ell^{\infty}\left(\Psi \right)$%
-valued martingale difference array with respect to filtration $%
\mathop{\mathcal F}\nolimits^T_{t}$, i.e. for every $t$, $\xi^T_{t}$ maps $%
\Omega^T$ to $\ell^{\infty}\left(\Psi \right)$, the space of bounded, $%
\mathop{\mathbb R}\nolimits$-valued functions on $\Psi$ with $\sup$-norm $\|
\cdot \|=\| \cdot \|_{\infty}$ and for each $u\in \Psi$, $\xi^T_{t}(u)$ is a 
$\mathop{\mathbb R}\nolimits$-valued martingale difference array: $%
\xi^T_{t}(u)$ is $\mathop{\mathcal F}\nolimits^T_{t}$-measurable and $%
\mathop{\rm E}\nolimits \left[\xi^T_{t}(u)\mid \mathop{\mathcal F}%
\nolimits^T_{t}\right]=0$. We are interested in studying the weak
convergence of discrete martingales $\sum_{t=1}^T\xi^T_{t}$. Denote a
decreasing series of finite partitions (DFP) of $\Psi$ as $\Pi=\left \{
\Pi(\varepsilon)\right \}_{\varepsilon \in (0,1) \cap\mathop{\mathbb Q}%
\nolimits}$, where $\Pi(\varepsilon)=\left \{ \Psi(\varepsilon;k)\right
\}_{1\le k \le N_{\Pi}(\varepsilon)}$ such that $\Psi=\bigcup_{k=1}^{N_{%
\Pi}(\varepsilon)} \Psi(\varepsilon;k)$, $N_{\Pi}(1)=1 $ and $%
\lim_{\varepsilon \to 0} N_{\Pi}(\varepsilon)=\infty$ monotonically in $%
\varepsilon$. The $\varepsilon$-entropy of the DFP $\Pi$ is $%
H_{\Pi}(\varepsilon)=\sqrt{\log N_{\Pi}(\varepsilon)}$. The quadratic $\Pi$%
-modulus of $\xi^T_{t}$ is $\mathop{\mathbb R}\nolimits_+\cup \{ \infty \}$%
-valued process 
\begin{equation}
\left \| \xi^T_{t}\right \|_{\Pi,T}=\sup_{\varepsilon \in (0,1) \cap%
\mathop{\mathbb Q}\nolimits}\frac{1}{\varepsilon}\max_{1\le k\le
N_{\Pi}(\varepsilon)}\sqrt{\sum_{t=1}^{T}\mathop{\rm E}\nolimits \left[%
\sup_{u,v\in \Psi(\varepsilon;k)}\left|\xi^T_{t}(u)-\xi^T_{t}(v)\right|^2
\mid\mathop{\mathcal F}\nolimits^T_{t}\right]}.
\end{equation}

\begin{theorema}
\label{thm:nish} Let $\{ \xi^T_{t}\}_{t=1,2,\ldots}$ be a $%
\ell^{\infty}\left(\Psi \right)$-valued martingale difference array and 
\newline
N1) (conditional variance convergence) $\sum_{t=1}^T\mathop{\rm E}\nolimits%
\left[\xi^T_{t}(u)\xi^T_{t}(v)\mid \mathop{\mathcal F}\nolimits^T_{t}\right]%
\to_{P^T}V(u,v)$ for every $u,v\in \Psi;$ \newline
N2) (Lindeberg condition) $\sum_{t=1}^T\mathop{\rm E}\nolimits \left[\left
\| \xi^T_{t}\right \|^2 1\left \{ \left \| \xi^T_{t}\right \|>{\varepsilon}%
\right \} \mid \mathop{\mathcal F}\nolimits^T_{t}\right]\to_{P^T}0$ for
every ${\varepsilon}>0$; \newline
N3) (partitioning entropy condition) there exist a DFP $\Pi$ of $\Psi$ such
that $\left \| \xi^T_{t}\right \|_{\Pi,T}=O_{P^T}(1)$ and $\int_0^1
H_{\Pi}(\varepsilon)d\varepsilon<\infty$. \newline
Then $\sum_{t=1}^T\xi^T_{t}\Rightarrow S$, where $S$ has normal marginals $%
\left(S\left(v_1\right),S\left(v_2\right),\ldots,S\left(v_d\right))\right)%
\sim_d N(0,\Sigma)$ with covariance $\Sigma=\left
\{V\left(v_i,v_j\right)\right \}_{ij}$.
\end{theorema}

\subsection{Additional technical assumptions}

To establish the asymptotic properties of the biparameter process $S_{2T}$
we need the following assumption for uniform convergence of different
empirical quantities.

\begin{assumptionA}
Under $H_{1T}$, the following uniform limits to continuous functions exist

\begin{enumerate}
\item $\mathop{\rm plim}_{T\rightarrow \infty }\frac{1}{T}%
\sum_{t=2}^{T}\gamma _{t-1,\theta _{0}}\left( u_{2},v_{2}\right) \gamma
_{t,\theta _{0}}\left( u_{1},v_{1}\right) $,

\item $\mathop{\rm plim}_{T\rightarrow \infty }\frac{1}{T}%
\sum_{t=2}^{T}I_{t-1,\theta _{0}}\left( v_{2}\right) \gamma _{t,\theta
_{0}}\left( u_{1},v_{1}\right) $,

\item $\mathop{\rm plim}_{T\rightarrow \infty }\frac{1}{T}%
\sum_{t=2}^{T}I_{t-1,\theta _{0}}\left( u_{2}\right) d\left( H_{t}\left(
\cdot \mid \Omega _{t}\right) ,F_{t,\theta _{0}}\left( \cdot \mid \Omega
_{t}\right) ,u_{1}\right) $, 

\item $\mathop{\rm plim}_{T\rightarrow \infty }\frac{1}{T}%
\sum_{t=2}^{T}I_{t-1,\theta _{0}}\left( u_{2}\right) \mathop{\rm E}\nolimits%
\left[ I_{t,\theta _{0}}\left( u_{1}\right) \ell _{t}\left( Y_{t},\Omega
_{t}\right) \mid \Omega _{t}\right] $,

\item $\mathop{\rm plim}_{T\rightarrow \infty }\frac{1}{T}%
\sum_{t=2}^{T}I_{t-1,\theta _{0}}\left( u_{2}\right) \nabla \left(
F_{t,\theta _{0}}\left( \cdot \mid \Omega _{t}\right) ,u_{1}\right) $.
\end{enumerate}
\end{assumptionA}

As it is discussed in the text, these conditions restrict the dynamics of
the data process such that some LLN holds, which is the case, e.g., for
stationary and ergodic processes.

\subsection{Proofs}

\begin{proof}[\textbf{Proof of Lemma~A}]
\textbf{(i)} By definition of $I_{F}\left( Y,u\right) $, $\mathop{\rm E}\nolimits_{G}%
\left[ I_{F}\left( Y,u\right) \right] =(1-\delta
_{F}(u))g(F^{-1}(u))+G(F^{-1}(u))-g(F^{-1}(u))=d\left( G,F,u\right) -\delta
_{F}(u)f(F^{-1}(u))+F(F^{-1}(u))=d\left( G,F,u\right) +u$. Similarly, by
direct calculation we obtain (ii), (iii), (vi) and (vii). We now provide a
detailed proof of (iv) and (v).

\textbf{(iv)} We prove a stronger result that for $G\in \mathop{\mathcal M}%
\nolimits$, such that $\sup_{k}\left \vert F\left( k\right) -G\left(
k\right) \right \vert \lor \left \vert H\left( k\right) -G\left( k\right)
\right \vert \leq \sup_{k}\left \vert F\left( k\right) -H\left( k\right)
\right \vert$ the expectation with respect to $G$ is bounded: $\mathop{\rm E}%
\nolimits_{G}\left[\left( I_{F}\left( Y,u\right) -I_{H}\left( Y,u\right)
\right)^2\right]\leq 9\sup_{k}\left \vert F\left( k\right) -H\left( k\right)
\right \vert$. Then, the required bound is obtained by setting $G\equiv F$.

Since $\left \vert I_{F}\left( Y,u\right) -I_{H}\left( Y,u\right) \right
\vert $ never exceeds $1 $, we have that $\mathop{\rm E}\nolimits_{G}\left[%
\left( I_{F}\left( Y,u\right) -I_{H}\left( Y,u\right) \right)^2\right]\leq %
\mathop{\rm E}\nolimits_{G}\left[\left \vert I_{F}\left( Y,u\right)
-I_{H}\left( Y,u\right) \right \vert \right]$, therefore we bound the latter
expectation.

Suppose that $F^{-1}(u)=H^{-1}(u)$. Then $I_{F}\left( Y,u\right)
-I_{H}\left( Y,u\right) =\delta _{H}(u)-\delta _{F}(u)$ for $Y=F^{-1}(u)$,
i.e. with probability $g\left( F^{-1}(u)\right) $, and is zero for other $Y$%
. Therefore, 
\begin{eqnarray*}
&&\mathop{\rm E}\nolimits_{G}\left[ \left \vert I_{F}\left( Y,u\right)
-I_{H}\left( Y,u\right) \right \vert \right] =\left \vert \delta
_{H}(u)-\delta _{F}(u)\right \vert g\left( F^{-1}(u)\right) \\
&\leq &\left \vert F\left( F^{-1}(u)\right) -H\left( F^{-1}(u)\right) \right
\vert +\left \vert f\left( F^{-1}(u)\right) -g\left( F^{-1}(u)\right) \right
\vert \delta _{F}(u)+\left \vert h\left( F^{-1}(u)\right) -g\left(
F^{-1}(u)\right) \right \vert \delta _{H}(u)\  \\
&\leq &\sup_{k}\left \vert F\left( k\right) -H\left( k\right) \right \vert
+\sup_{k}\left \vert f\left( k\right) -g\left( k\right) \right \vert
+\sup_{k}\left \vert h\left( k\right) -g\left( k\right) \right \vert \\
&\leq &5\sup_{k}\left \vert F\left( k\right) -H\left( k\right) \right \vert
,\ 
\end{eqnarray*}%
since $\delta _{F}(u),\delta _{H}(u)\in \lbrack 0,1)$ and $\sup_{k}\left
\vert h\left( k\right) -g\left( k\right) \right \vert \leq 2\sup_{k}\left
\vert F\left( k\right) -G\left( k\right) \right \vert$.

Suppose that $F^{-1}(u)<H^{-1}(u)$. Note that $I_{F}\left( Y,u\right)
-I_{H}\left( Y,u\right) =0$ for $Y\not \in \lbrack F^{-1}(u),H^{-1}(u)]$. We
separately bound each term in 
\begin{eqnarray*}
\mathop{\rm E}\nolimits_{G}\left[ \left \vert I_{F}\left( Y,u\right)
-I_{H}\left( Y,u\right) \right \vert \right] &=&\mathop{\rm E}\nolimits_{G}%
\left[ \left \vert I_{F}\left( Y,u\right) -I_{H}\left( Y,u\right) \right
\vert \mathbbm{1}\! \left \{ Y=F^{-1}(u)\right \} \right] \\
&&+\mathop{\rm E}\nolimits_{G}\left[ \left \vert I_{F}\left( Y,u\right)
-I_{H}\left( Y,u\right) \right \vert \mathbbm{1}\! \left \{
Y=H^{-1}(u)\right \} \right] \\
&&+\mathop{\rm E}\nolimits_{G}\left[ \left \vert I_{F}\left( Y,u\right)
-I_{H}\left( Y,u\right) \right \vert \mathbbm{1}\! \left \{
F^{-1}(u)<Y<H^{-1}(u)\right \} \right] .
\end{eqnarray*}

For $Y=F^{-1}(u)$, $I_{F}\left( Y,u\right) -I_{H}\left( Y,u\right) =-\delta
_{F}(u)$. Then 
\begin{eqnarray*}
&&\mathop{\rm E}\nolimits_{G}\left[ \left( I_{F}\left( Y,u\right)
-I_{H}\left( Y,u\right) \right) \mathbbm{1}\! \left \{ Y=F^{-1}(u)\right \} %
\right] =\left \vert I_{F}\left( F^{-1}(u),u\right) -I_{H}\left(
F^{-1}(u),u\right) \right \vert g(F^{-1}(u)) \\
&=&\delta _{F}(u)g\left( F^{-1}(u)\right) \\
&=&F\left( F^{-1}(u)\right) -u+\delta _{F}(u)\left( g\left( F^{-1}(u)\right)
-f\left( F^{-1}(u)\right) \right) \\
&\leq &\sup_{k}\left \vert F\left( k\right) -H\left( k\right) \right \vert
+\sup_{k}\left \vert f\left( k\right) -g\left( k\right) \right \vert \leq
3\sup_{k}\left \vert F\left( k\right) -H\left( k\right) \right \vert ,
\end{eqnarray*}
since $\delta _{F}(u)\in \lbrack 0,1)$ and for $u\in \lbrack H\left(
F^{-1}(u)\right) ,F\left( F^{-1}(u)\right) ]$ we have that $F\left(
F^{-1}(u)\right) -u\leq F\left( F^{-1}(u)\right) -H\left( F^{-1}(u)\right) $.

For $Y=H^{-1}(u)$, $I_{F}\left( Y,u\right) -I_{H}\left(
Y,u\right)=-1+\delta_H(u)$.

Then 
\begin{eqnarray*}
&&\mathop{\rm E}\nolimits_{G}\left[ \left \vert I_{F}\left( Y,u\right)
-I_{H}\left( Y,u\right) \right \vert \mathbbm{1}\! \left \{
Y=H^{-1}(u)\right \} \right] \\
&=&\left \vert I_{F}\left( H^{-1}(u),u\right) -I_{H}\left(
H^{-1}(u),u\right) \right \vert g(H^{-1}(u)) \\
&=&\left( 1-\delta _{H}(u)\right) g\left( H^{-1}(u)\right) \\
&=&u-H\left( H^{-1}(u)-1\right) +\left( 1-\delta _{H}(u)\right) \left(
g\left( H^{-1}(u)\right) -h\left( H^{-1}(u)\right) \right) \\
&\leq &\sup_{k}\left \vert F\left( k\right) -H\left( k\right) \right \vert
+\sup_{k}\left \vert h\left( k\right) -g\left( k\right) \right \vert \leq
3\sup_{k}\left \vert F\left( k\right) -H\left( k\right) \right \vert ,
\end{eqnarray*}
since $\delta _{H}(u)\in \lbrack 0,1)$ and for $u\in \lbrack H\left(
H^{-1}(u)-1\right) ,F\left( H^{-1}(u)-1\right) ]$ we have that $u-H\left(
H^{-1}(u)-1\right) \leq F\left( H^{-1}(u)-1\right) -H\left(
H^{-1}(u)-1\right) $.

For $F^{-1}(u)<Y<H^{-1}(u)$, $I_{F}\left( Y,u\right) -I_{H}\left( Y,u\right)
=-1$. Then 
\begin{eqnarray*}
&&E_{F}\left \vert I_{F}\left( Y,u\right) -I_{H}\left( Y,u\right) \right
\vert \mathbbm{1}\! \left \{ F^{-1}(u)<Y<H^{-1}(u)\right \} \\
&=&\sum_{k=F^{-1}(u)+1}^{H^{-1}(u)-1}g(k)=G(H^{-1}(u)-1)-G(F^{-1}(u)) \\
&\leq &F(H^{-1}(u)-1)-F(F^{-1}(u))+2\sup_{k}\left \vert G(k)-F(k)\right \vert
\\
&\leq &F(H^{-1}(u)-1)-H(H^{-1}(u)-1)+2\sup_{k}\left \vert G(k)-F(k)\right
\vert \leq 3\sup_{k}\left \vert F(k)-H(k)\right \vert
\end{eqnarray*}
since $H(H^{-1}(u)-1)<u<F(F^{-1}(u))<F(H^{-1}(u)-1)$.

Adding everything together, get that $\mathop{\rm E}\nolimits_{G}\left[\left
\vert I_{F}\left( Y,u\right) -I_{H}\left( Y,u\right) \right \vert^2\right]
\le9\sup_{k}{\left \vert {F}\left( k\right) -{H}\left( k\right) \right \vert 
}$ for $F^{-1}(u)<H^{-1}(u)$. This equation is symmetric with respect to $F$
and $H$; therefore, it holds also for $F^{-1}(u)>H^{-1}(u)$.

\textbf{(v)} Let $[a,b]$ denote the interval $\Psi (\varepsilon )$ of length 
$\varepsilon ^{2}$, $\sup \xi ^{2}$ denote the supremum of $\xi ^{2}$ over $%
u,v\in \lbrack a,b]$, where $\xi :=I_{F}\left( Y,u\right) -u-I_{F}\left(
Y,v\right) +v$.

Note that $\left \vert \xi \right \vert \leq 1$; moreover, if $%
[F(Y-1),F(Y)]\cap \lbrack a,b]=\emptyset $, then $\sup \left \vert \xi
\right \vert =\varepsilon ^{2}$ and if $[a,b]\subset \lbrack F(Y-1),F(Y)]$,
then $\sup \left \vert \xi \right \vert =\frac{1-f(Y)}{f(Y)}\varepsilon ^{2}$%
.

Suppose that $F^{-1}(a)=F^{-1}(b)$, i.e. $[a,b]\subset \lbrack
F(F^{-1}(a)-1),F(F^{-1}(a))]$. Then $\mathop{\rm E}\nolimits_{F}\left[ \sup
\xi ^{2}\right] \leq \mathop{\rm E}\nolimits_{F}\left[ \sup \left \vert \xi
\right \vert \right] =\varepsilon ^{2}\sum_{k\neq F^{-1}(a)}f(k)+\left( 
\frac{1-f(F^{-1}(a))}{f(F^{-1}(a))}\right) \varepsilon
^{2}f(F^{-1}(a))=2(1-f(F^{-1}(a)))\varepsilon ^{2}\leq 2\varepsilon ^{2}$.

Suppose that $F^{-1}(a)<F^{-1}(b)$, i.e. $[a,b]$ contains at least one point 
$F(k)$ or even intervals $[F(k-1),F(k)]\subset[a,b]$. On such intervals, $%
|\xi|$ goes up to $1-f(k)$, but the probability of $Y$ taking all such $k$
is bounded by $b-a$. More precisely, 
\begin{eqnarray*}
\mathop{\rm E}\nolimits_{F}\left[ \sup \xi ^{2}\right] &\leq &\mathop{\rm E}%
\nolimits_{F}\left[ \sup \left \vert \xi \right \vert \right] \\
&=&\varepsilon ^{2}\sum_{k<F^{-1}(a)-1}f(k)+\left( \frac{1-f(F^{-1}(a))}{%
f(F^{-1}(a))}\right) \varepsilon ^{2}f(F^{-1}(a)) \\
&&+\sum_{k\in \lbrack F^{-1}(a)+1,F^{-1}(b)-1]}f(k)+\left( \frac{%
1-f(F^{-1}(b))}{f(F^{-1}(b))}\right) \varepsilon
^{2}f(F^{-1}(b))+\varepsilon ^{2}\sum_{k>F^{-1}(b)}f(k) \\
&<&4\varepsilon ^{2},
\end{eqnarray*}
since the sum of the first and the last terms is below $\varepsilon ^{2}$,
the second and the fourth terms each is bounded by $\varepsilon ^{2}$ and
the third term is $\sum_{k\in \lbrack F^{-1}(a)+1,F^{-1}(b)-1]}f(k)=F\left(
F^{-1}(b)-1\right) -F\left( F^{-1}(a)+1\right) \leq b-a=\varepsilon ^{2}$.
\end{proof}

\begin{proof}[\textbf{Proof of Lemma~\protect\ref{lem:mds}}]
Substitute $G=F=F_{\theta _{0}}\left( \cdot \mid \Omega _{t}\right) $ in
Lemma~A(\ref{a:ei}) to demonstrate that $E\left[ I_{t,\theta _{0}}\left(
u\right) \mid \Omega _{t}\right] =E\left[ I_{t,\theta _{0}}\left( u\right) %
\right] =u$, therefore $I_{t,\theta _{0}}\left( u\right) -u$ is a martingale
difference sequence for every $u\in \left[ 0,1\right] $. The conditional
variance expression follows from Lemma~A(\ref{a:eii}) by taking $%
G=F=F_{\theta _{0}}\left( \cdot \mid \Omega _{t}\right) $. 
\end{proof}

However the $I_{t,\theta _{0}}\left( u\right) $ are not independent in
general. To show that, note that bivariate independence requires that%
\begin{equation*}
\Pr \left( I_{t,\theta _{0}}\left( u\right) \leq u_{1},I_{t-1,\theta
_{0}}\left( u\right) \leq u_{2}\right) =\Pr \left( I_{t,\theta _{0}}\left(
u\right) \leq u_{1}\right) \Pr \left( I_{t-1,\theta _{0}}\left( u\right)
\leq u_{2}\right) 
\end{equation*}%
for all $u,$ $u_{1}$ and $u_{2}\in \left[ 0,1\right] $. Now we see that the
lhs is 
\begin{equation*}
\begin{split}
\mathop{\rm E}\nolimits\left[ \mathbbm{1}\!\left\{ I_{t,\theta _{0}}\left(
u\right) \leq u_{1}\right\} \mathbbm{1}\!\left\{ I_{t-1,\theta _{0}}\left(
u\right) \leq u_{2}\right\} \right] ={}& \mathop{\rm E}\nolimits\left[ %
\mathop{\rm E}\nolimits\left[ \mathbbm{1}\!\left\{ I_{t,\theta _{0}}\left(
u\right) \leq u_{1}\right\} \mathbbm{1}\!\left\{ I_{t-1,\theta _{0}}\left(
u\right) \leq u_{2}\right\} \mid \Omega _{t}\right] \right]  \\
={}& \mathop{\rm E}\nolimits\left[ \mathbbm{1}\!\left\{ I_{t-1,\theta
_{0}}\left( u\right) \leq u_{2}\right\} \mathop{\rm E}\nolimits\left[ %
\mathbbm{1}\!\left\{ I_{t,\theta _{0}}\left( u\right) \leq u_{1}\right\}
\mid \Omega _{t}\right] \right] 
\end{split}%
\end{equation*}%
and now, for $u_{1},u\in \left( 0,1\right) $ and under $H_{0},$ 
\begin{equation*}
\begin{split}
\mathop{\rm E}\nolimits\left[ \mathbbm{1}\!\left\{ I_{t,\theta _{0}}\left(
u\right) \leq u_{1}\right\} \mid \Omega _{t}\right] ={}& 1-F_{\theta
_{0}}\left( F_{\theta _{0}}^{-1}\left( u\mid \Omega _{t}\right) \mid \Omega
_{t}\right)  \\
& +\mathbbm{1}\!\left\{ 1-\delta _{F_{\theta _{0}}\left( \cdot \mid \Omega
_{t}\right) }\left( u\right) \leq u_{1}\right\} f{_{\theta _{0}}\left( \cdot
\mid \Omega _{t}\right) }\left( F_{\theta _{0}}^{-1}\left( u\mid \Omega
_{t}\right) \right) ,
\end{split}%
\end{equation*}%
which depends on $\Omega _{t},$ and therefore $E\left( \mathbbm{1}\!\left\{
I_{t,\theta _{0}}\left( u\right) \leq u_{1}\right\} \mid \Omega _{t}\right)
\neq E\left( \mathbbm{1}\!\left\{ I_{t,\theta _{0}}\left( u\right) \leq
u_{1}\right\} \right) $ with positive probability, and independence does not
follow in general.

\begin{proof}[\textbf{Proof of Lemma \protect\ref{lem:comp}.}]
Because $U_{t}^{r}\left( {\theta }_{0}\right) $ are continuous, $\hat{F}_{{%
\theta }_{0}}^{r}\left( u\right) $ is a (uniform) consistent estimate of cdf
of $U_{t}^{r}\left( {\theta }_{0}\right) $. Then by Lemma~A(\ref{a:eid}) and
A(\ref{a:eiim}) and ULLN we get the uniform consistency of $\hat{F}_{{\theta 
}_{0},M}^{r}\left( u\right) $ and $\tilde{F}_{{\theta }_{0}}^{r}\left(
u\right) $. The efficiency gain comes from Lemma~A(\ref{a:ii}).
\end{proof}

\begin{proof}[\textbf{Proof of Lemma \protect\ref{lem:limitS}.}]
We need to verify conditions N1-N3 of Theorem \ref{thm:nish}. Fix $%
\varepsilon >0$ and take $\Psi =[0,1]$ with usual norm and equidistant
partition $0=u_{0}<u_{1}<\ldots <u_{N_{\Pi }\left( \varepsilon \right) }=1$,
i.e. partition of $[0,1]$ in $N_{\Pi }\left( \varepsilon \right)
=[\varepsilon ^{-2}]+1$ equal intervals of length $\varepsilon ^{2}$ (the
last interval maybe even smaller), $\Psi (\varepsilon ;k)=[u_{k-1},u_{k}]$
and $\xi _{t}^{T}=\left( I_{F}\left( Y_{t},u\right) -u\right) /\sqrt{T}$,
which is a square integrable martingale difference by Lemma \ref{lem:mds}.
Then Condition N1 follows from Lemma \ref{lem:mds} and Assumption 1.
Condition N2 is satisfied because for $T>1+\left[ {\varepsilon }^{-2}\right] 
$, the indicator $1\left \{ \sup_{u\in \lbrack 0,1]}\left \vert I_{F}\left(
Y_{t},u\right) -u\right \vert /\sqrt{T}>{\varepsilon }\right \} =0$.
Condition N3 follows from the bound in Lemma~A(\ref{a:i-i}). Indeed, $%
\int_{0}^{1}H_{\Pi }(\varepsilon )d\varepsilon <\infty $ and 
\begin{equation*}
\left \Vert \xi _{t}^{T}\right \Vert _{\Pi ,k}\leq \sup_{\varepsilon \in
(0,1)\cap \mathop{\mathbb Q}\nolimits}\frac{1}{\varepsilon }\max_{1\leq
k\leq N_{\Pi }(\varepsilon )}\sqrt{\varepsilon ^{2}}\leq 1\quad \text{a.s.}
\end{equation*}%

\end{proof}

\begin{proof}[\textbf{Proof of Lemma \protect\ref{lem:limitSLA}.}]
Apply weak convergence result from Lemma \ref{lem:limitS} under $G_{T,\theta
_{0}}\left( \cdot \mid \Omega _{t}\right)$ with $\xi^T_t:=\left(I_{F_{\theta
_{0}}\left( \cdot \mid \Omega
_{t}\right)}\left(Y_t,u\right)-u-d\left(G_{T,\theta _{0}}\left( \cdot \mid
\Omega _{t}\right),F_{\theta _{0}}\left( \cdot \mid \Omega
_{t}\right),u\right)\right)/\sqrt{T}$, which is a \newline
square integrable martingale difference because of Lemma~A(\ref{a:ei}) with $%
G=G_{T,\theta _{0}}\left( \cdot \mid \Omega _{t}\right)$ and $F=F_{\theta
_{0}}\left( \cdot \mid \Omega _{t}\right)$. Then Condition N1 follows from
Lemma~A(\ref{a:eii}) and the fact that $d\left(G,F, u,v\right)$ are bounded
in absolute value by $T^{-1/2}$ a.s. Condition N2 is satisfied because for $%
T>1+\left[{\varepsilon}^{-2}\right]$, the indicator is $0$. Condition N3
follows from the bound in Lemma~A(\ref{a:i-i}) and the fact that $\left(%
\mathop{\rm E}\nolimits_{G}\left[\cdot \right]-\mathop{\rm E}\nolimits_{F}%
\left[\cdot \right]\right)$ applied to a.s. bounded r.v. are bounded in
absolute value by $T^{-1/2}$ a.s. We obtain that $\sum_{t=1}^T\xi^T_{t}%
\Rightarrow S$, the same limit as in Lemma \ref{lem:limitS}. Finally, use
additivity of $d\left(\cdot,\cdot,\cdot \right)$ in the first argument and
apply ULLN to $S_T-\sum_{t=1}^T\xi^T_{t}=\sum_{t=1}^T d\left(G_{T,\theta
_{0}}\left( \cdot \mid \Omega _{t}\right),F_{\theta _{0}}\left( \cdot \mid
\Omega _{t}\right),u\right)/\sqrt{T}=\delta \sum_{t=1}^T d\left(H\left(
\cdot \mid \Omega _{t}\right),F_{\theta _{0}}\left( \cdot \mid \Omega
_{t}\right),u\right)/{T}$.
\end{proof}

\begin{proof}[\textbf{Proof of Lemma \protect\ref{lem:puexpansion}.}]
Under $H_{1T}$, i.e. under $G_{T,\theta_0}$, Equation (\ref{eq:puexpansion})
can be established using standard methods, applying Doob and Rosenthal
inequalities for MDS (Hall and Heyde, 1980) 
$\sqrt{T}\xi^T_t:=
I_{F_{\hat \theta_T} \left( \cdot \mid \Omega _{t}\right)}\left(Y_t,u\right)
-I_{F_{\theta _{0}}\left( \cdot \mid \Omega _{t}\right)}\left(Y_t,u\right)
-d\left(G_{T,\theta_0}\left( \cdot \mid \Omega _{t}\right),F_{\hat
\theta_T}\left( \cdot \mid \Omega _{t}\right),u\right)$ $+d\left(G_{T,%
\theta_0}\left( \cdot \mid \Omega _{t}\right),F_{\theta_0}\left( \cdot \mid
\Omega _{t}\right),u\right).
$ 
Define $z_T:=\sum_{t=1}^T \xi^T_t$. When it is necessary, we will write
explicitly arguments: $z_{T}(u,\hat \theta_T)$. We show that $%
\sup_u\left|z_T\right|=o_p(1)$. Since\newline
$\sqrt{T}\left(\hat \theta_T-\theta_0\right)=O_P(1)$, it is sufficient to
establish that for some $\gamma<1/2$ 
\begin{equation*}  \label{eq:ngamma}
\sup_{u,\left \Vert \eta -\theta _{0}\right \Vert \leq T^{-\gamma }}\left
\vert z_{T}(u,\eta )\right \vert =o_{p}(1).
\end{equation*}%
Note that for $T>\delta^2/\nu_1^2$, by Assumption 3C, 
\begin{eqnarray}  \label{eq:G-Fhat}
\Pr \left(\sup_{\eta,t}\max_{y}\left|G_{T,t,\theta_0}\left(y\mid \Omega
_{t}\right)-F_{t,\eta}\left( y\mid \Omega _{t}\right)\right|>\nu_1\right)\le
M_F T^{-\gamma}/\nu_1.
\end{eqnarray}

First, we will show that $\forall \ \eta, u\
\left|z_{T}\right|=o_p\left(1\right)$. Since $\xi^T_{t}$ are bounded by 2 in
absolute value and form a martingale difference sequence with respect to $%
\Omega _{t}$, by the Doob inequality $\forall p\ge 1$ and $\forall
\varepsilon>0$%
\begin{equation*}
P\left( \max_{t=1, \ldots,T}\left \vert z_{t} \right \vert >\varepsilon
\right) \leq E \left| z_{T}\right| ^{p}/\varepsilon ^{p},
\end{equation*}%
and by Rosenthal 
inequality, $\forall p\ge 2\ \exists C_1$%
\begin{equation*}
E \left| z_{T}\right| ^{p} \leq C_1\left[ E\left \{ \sum E\left(
\left(\xi_t^{T}\right)^{2} \mid \Omega _{t}\right) \right \} ^{p/2} +\sum
E\left| \xi^T_t \right| ^{p}\right].
\end{equation*}%
Take $p=4$. The first term is small because of bounds in Lemma~A(\ref{a:ai-i}%
) and (\ref{eq:G-Fhat}). Because $\left|\xi^T_t\right|\le 2/\sqrt{T}$, $\sum
E\left| \xi^T_t \right| ^{p}\le 2T^{1-p/2}$. Therefore we have a pointwise
bound. Uniformity in $u,\eta$ can be established using monotonicity of $%
I_{F_{\theta}\left( \cdot \mid \Omega _{t}\right)}\left(Y_t,u\right)$ and
continuity of $d\left(G_{T,\theta_0}\left( \cdot \mid \Omega
_{t}\right),F_{\hat \theta_T}\left( \cdot \mid \Omega _{t}\right),u\right)$
by employing bounds in Lemma~A(\ref{a:ai-i}) and (\ref{eq:G-Fhat}).

Finally, use that uniformly in $u$ 
\begin{equation*}
\begin{split}
\frac{1}{\sqrt T}\sum_t\left( d\left(G_{T,\theta_0}\left( \cdot \mid \Omega
_{t}\right),F_{\hat \theta_T}\left( \cdot \mid \Omega _{t}\right),u\right)
-d\left(G_{T,\theta_0}\left( \cdot \mid \Omega
_{t}\right),F_{\theta_0}\left( \cdot \mid \Omega _{t}\right),u\right)\right)
\\
=\sqrt{T}\left(\hat \theta_T-\theta_0\right) \frac{1}{T}\sum_t \nabla
\left(F_{\theta_0}\left( \cdot \mid \Omega _{t}\right),u\right) +o_p(1).
\end{split}%
\end{equation*}
\end{proof}


\begin{proof}[\textbf{Proof of Theorem \protect\ref{thm:limitNormSLAPEE}.}]
The joint weak convergence (\ref{eq:pujoint}) follows from
finite-dimensional convergence by CLT for MDS, while tightness was
established in the proof of Lemma~\ref{lem:limitSLA}.
\end{proof}


\begin{proof}[\textbf{Proof of Theorem \protect\ref{thm:limitNormSLAPEE2}.}]

Note that 
\begin{equation*}
S_{2T}=\sum_{t=2}^{T}\xi _{t}^{T}+\frac{1}{(T-1)^{1/2}}\left\{ \left(
I_{\theta _{0},T}\left( u_{1}\right) -u_{1}\right) I_{\theta _{0},T-1}\left(
u_{2}\right) +u_{1}\left( I_{\theta _{0},1}\left( u_{2}\right) -u_{2}\right)
\right\} ,
\end{equation*}%
where 
\begin{equation*}
\xi _{t}^{T}:=\frac{1}{(T-1)^{1/2}}\left\{ \left( I_{t,\theta _{0}}\left(
u_{1}\right) -u_{1}\right) I_{t-1,\theta _{0}}\left( u_{2}\right)
+u_{1}\left( I_{t,\theta _{0}}\left( u_{2}\right) -u_{2}\right) \right\} 
\end{equation*}%
is a square integrable martingale difference by Lemma \ref{lem:mds}. The
rest is similar to the proof of Theorem~\ref{thm:limitNormSLAPEE}. To obtain 
$S_{2T}\left( u\right) \Rightarrow S_{2\infty }\left( u\right) $ under $H_{0}
$, verify conditions N1-N3 of Theorem~\ref{thm:nish} for $\xi _{t}^{T}$ as
it is done in the proof of Lemma \ref{lem:limitS}. The covariance function
of $S_{2\infty }\left( u\right) $ is 
\begin{equation*}
\begin{split}
V_{2}\left( u,v\right) :={}& \left( u_{1}\wedge v_{1}\right) \left(
u_{2}\wedge v_{2}\right) -3u_{1}v_{1}u_{2}v_{2} \\
& +\left( u_{1}\wedge v_{1}\right) \mathop{\rm plim}_{T\rightarrow \infty }%
\frac{1}{T}\sum_{t=2}^{T}\gamma _{t-1,\theta _{0}}\left( u_{2},v_{2}\right) 
\\
& -\mathop{\rm plim}_{T\rightarrow \infty }\frac{1}{T}\sum_{t=2}^{T}\gamma
_{t,\theta _{0}}\left( u_{1},v_{1}\right) \left( I_{t-1,\theta _{0}}\left(
u_{2}\wedge v_{2}\right) -\gamma _{t-1,\theta _{0}}\left( u_{2},v_{2}\right)
\right)  \\
& +\left( u_{2}\wedge v_{1}\right) u_{1}v_{2}-u_{1}\mathop{\rm plim}%
_{T\rightarrow \infty }\frac{1}{T}\sum_{t=2}^{T}\gamma _{t,\theta
_{0}}\left( u_{1},v_{1}\right) I_{t-1,\theta _{0}}\left( v_{2}\right)  \\
& +\left( u_{1}\wedge v_{2}\right) u_{2}v_{1}-v_{1}\mathop{\rm plim}%
_{T\rightarrow \infty }\frac{1}{T}\sum_{t=2}^{T}\gamma _{t,\theta
_{0}}\left( u_{1},v_{1}\right) I_{t-1,\theta _{0}}\left( u_{2}\right) .
\end{split}%
\end{equation*}

Under $H_{1T}$, apply the same weak convergence result under $G_{T,t,\theta
_{0}}\left( \cdot \mid \Omega _{t}\right) $ with 
\begin{equation*}
\begin{split}
\zeta _{t}^{T}:={}& \xi _{t}^{T}-I_{t-1,\theta _{0}}\left( u_{2}\right)
d\left( G_{T,t,\theta _{0}}\left( \cdot \mid \Omega _{t}\right) ,F_{t,\theta
_{0}}\left( \cdot \mid \Omega _{t}\right) ,u_{1}\right) /\sqrt{T-1} \\
& +u_{1}d\left( G_{T,t,\theta _{0}}\left( \cdot \mid \Omega _{t}\right)
,F_{t,\theta _{0}}\left( \cdot \mid \Omega _{t}\right) ,u_{2}\right) /\sqrt{%
T-1},
\end{split}%
\end{equation*}%
which is a square integrable martingale difference because of Lemma~A(\ref%
{a:ei}) with $G=G_{T,t,\theta _{0}}\left( \cdot \mid \Omega _{t}\right) $
and $F=F_{\theta _{0}}\left( \cdot \mid \Omega _{t}\right) $. Then proceed
as in the proof of Lemma \ref{lem:limitSLA}.

In order to establish (\ref{eq:puexpansion2}), repeat the steps of the proof
of Lemma \ref{lem:puexpansion} for $\tilde
\zeta^{T}_t:=\zeta^{T}_t-\hat\zeta^{T}_t$, where $\hat \zeta^{T}_t$ is $%
\zeta^{T}_t$ with $F_{t,\hat \theta_T}$ in place of $F_{t,\theta_0}$.
\end{proof}

\begin{proof}[\textbf{Proof of Theorem \protect\ref{thm:limitNormSBoot}.}]
Repeat the arguments of the proofs of Theorems \ref{thm:limitNormSLAPEE} and %
\ref{thm:limitNormSLAPEE2} for sample generated by $F_{\theta _{T}}$,
defined in Assumption~6, to obtain conditional convergence. Then follow as
in Andrews (1997) proof of Corollary~1. 
\end{proof}

\subsection{Checking assumptions for the Poisson model}

Here we write $Y_t$ for $Y_{t}^{\star}$. For Poisson model $Y_{t}\mid \Omega _{t}\sim \text{Poisson}(\lambda _{t})$ the probability distribution is $\Pr(Y_t=k\mid\Omega_t)=P_{\lambda_t}(k)=\frac{\lambda_t^{k}\exp(-\lambda_t)}{k!}$
 and the cumulative distribution function is 
\begin{equation*}
F_{t,\theta}(k\mid\Omega_t)=\sum_{j=0}^k\Pr(Y_t=j\mid\Omega_t)=\sum_{j=0}^k\frac{\lambda_t^{j}\exp(-\lambda_t)}{j!}=Q(k,\lambda_t),
\end{equation*}
 where $Q(\cdot,\cdot)$ is the regularized gamma function, and $\lambda_t=\lambda_t(\beta)=\exp(X_t'\beta)$, $t=1,2,\ldots$. If covariates $X_t$ are iid or stationary and ergodic, and $\Omega_t$ omits lags of the dependent variable $Y_t,$ then the LLN applies both under the null and local alternatives (like, e.g., the local alternative considered in Eq. (2.12) in Cameron and Trivedi, 1990) to justify Assumptions 2-6 and Assumption A, which involve functions of $\Omega_t$ that are uniformly continuous in $u$. However, it can also be interesting to allow the intensity to depend on lags of the dependent variable. For simplicity we consider $AR(1)$ dynamics. $AR(p)$ can be treated similarly but is more lengthy. The parameters enter through $\lambda_t=\lambda_t(\theta)=\alpha_0+\alpha_1\lambda_{t-1}+\rho Y_{t-1}$, $t=1,2,\ldots$, and are gathered in $\theta=(\alpha_0, \alpha_1, \rho)'$. We assume that $\alpha_0, \alpha_1, \rho$ are positive, $\lambda_0$ and $Y_0$ are fixed and $\alpha_1+\rho<1$. Under these conditions, there exist a unique stationary and ergodic solution to this model (Fokianos et al., 2009). Such data generating processes allow to use results on (generic, uniform) LLN, which facilitate the checking of assumptions in the paper.  Conditions for stationarity and ergodicity for nonlinear $\lambda_t(\theta)$ can be found in Neumann (2011) and are directly applicable to the analysis under the null hypothesis. 
However, we are not aware of LLN results for these models under local alternatives despite Fokianos and Neumann (2013, Proposition 2.3(ii)) use related arguments.

Let $\lambda_{t,0}=\lambda_t(\theta_0)$ and the null hypothesis is $Y_{t}\mid \Omega _{t}\sim \text{Poisson}(\lambda _{t,0})$ for some $\theta_0\in\Theta$. Then $U_t=Q(Y_t+1,\lambda_{t,0})$ and $U^{-}_t=Q(Y_t,\lambda_{t,0})$, and the \emph{nonrandomized transform} $%
Y_{t}\mapsto I_{t,\theta _{0}}\left( u\right) $ for $u\in \lbrack 0,1]$ is 
\begin{equation*}
I_{t,\theta _{0}}\left( u\right) =\left\{ 
\begin{array}{rrr}
0, &  & u\leq Q(Y_t,\lambda_{t,0}); \\ 
\displaystyle{\ \frac{u-Q(Y_t,\lambda_{t,0})}{\lambda_{t,0}^{Y_t}\exp(-\lambda_{t,0})}Y_t !,} &  & 
Q(Y_t,\lambda_{t,0}) \leq u\leq Q(Y_t+1,\lambda_{t,0})
; \\ 
1, &  & Q(Y_t+1,\lambda_{t,0})\leq u,%
\end{array}%
\right.   
\end{equation*}
from where one obtains the empirical processes and the test statistics defined in Sections~1-2.

Now consider Assumption 1. For Poisson model
\begin{equation*}
\gamma _{t,\theta _{0}}\left( u,v\right) =\frac{\left( Q(k+1,\lambda_{t,0})-u\vee v\right)
\left( u\wedge v-Q(k,\lambda_{t,0})\right) }{\lambda_{t,0}^{k}\exp(-\lambda_{t,0})}k !\mathbbm{1}\!\left\{ k(u)=k(v)\right\},
\end{equation*}%
where $k=k\left( u\right) =\min \{y:Q(y,\lambda_{t,0})\geq u\}$. For the Poisson DGP described above, $Y_t$ is stationary and ergodic, $\gamma _{\infty
}\left( u,v\right) :=\mathop{\rm E}\nolimits \left[ \gamma _{1,\theta
_{0}}\left( u,v\right) \right]$ satisfies Assumption 1. By the same argument Assumptions 2, 3D, 4C, 5 are fulfilled.

Assumption 3A and 3B are trivial. For Assumption 3C note that 
\begin{equation*}
\dot{F}_{t,\theta }\left( k\mid \Omega _{t}\right) = \left(\sum_{j=0}^{k-1}\frac{\lambda_t^{j}}{j!}-\sum_{j=0}^k\frac{\lambda_t^{j}}{j!}\right)	\exp(-\lambda_t)\dot{\lambda_{t}}=-\frac{\lambda_t^{k}}{k!}\exp(-\lambda_t)\dot{\lambda_{t}},
\end{equation*}
where
\begin{equation*}
\dot{\lambda_{t}}=\left(1+\alpha_1 \frac{\partial\lambda_{t-1}}{\partial\alpha_0},\lambda_{t-1}+\alpha_1\frac{\partial\lambda_{t-1}}{\partial\alpha_1},Y_{t-1}+\alpha_1\frac{\partial\lambda_{t-1}}{\partial\rho}\right)'.
\end{equation*}
The last expression can be iterated from $t-1$ to $t=1$ and because $\alpha_1<1$ the arithmetic progression sum of mean squares is bounded, as in the proof of Lemma 3.2 of Fokianos et al. (2009).

Assumption 4A, 4B and 6B are standard, see e.g. Andrews (1997) which adapts to Poisson model using Theorem 3.1 of Fokianos et al. (2009).

Assumption 6A is trivial, because there is no explanatory variables other than own past values.

\section{Acknowledgements}
We thank Juan Mora for useful comments. 
Support from the Ministerio Economia y Competitividad (Spain), grants ECO2012-31748, ECO2014-57007p, MDM 2014-0431, Comunidad de Madrid, MadEco-CM (S2015/HUM-3444), and  Fundaci\'{o}n Ram\'{o}n Areces is gratefully acknowledged. 

\renewcommand{\baselinestretch}{1}

\end{document}